\documentstyle[amssymb,amsfonts]{amsart}

\newenvironment{proof}{\noindent {\bf Proof} }{\endprf\par}
\def \endprf{\hfill  {\vrule height6pt width6pt depth0pt}\medskip}
\def\emph#1{{\it #1}}
\def\textbf#1{{\bf #1}}

\newcommand{\bea}{\begin{eqnarray}}
\newcommand{\eea}{\end{eqnarray}}
\def\beaa{\begin{eqnarray*}}
\def\eeaa{\end{eqnarray*}}

\def\ba{\begin{array}}
\def\ea{\end{array}}
\def\be#1{\begin{equation} \label{#1}}
\def \eeq{\end{equation}}

\newcommand{\nn}{\nonumber}

\def\nn{\nonumber}

\def\D{{\cal D}}

\def\MM{{\cal M}}
\def\Ml{{\cal M}_{\ell}}

\def\a{\alpha}
\def\b{\beta}
\def\ga{\gamma}

\def\Ga{\Gamma}
\def\de{\delta}

\def\la{\lambda}
\def\La{\Lambda}
\def\si{\sigma}

\def\nab{\nabla}

\def\f14{\frac{1}{4}}
\def\f12{{\frac{1}{2}}}

\def\c{\cdot}

\newcommand{\les}{\lesssim}

\newcommand{\DD}{{\mathcal D}}

\def\lap{\Delta}
\def\div{\mbox{ div }}
\def\curl{\mbox{ curl }}
\def\dual{{\,^*}}

\def\pr{\partial}

\begin{document}
\theoremstyle{plain}
  \newtheorem{theorem}[subsection]{Theorem}
  \newtheorem{conjecture}[subsection]{Conjecture}
  \newtheorem{proposition}[subsection]{Proposition}
  \newtheorem{lemma}[subsection]{Lemma}
  \newtheorem{corollary}[subsection]{Corollary}

\theoremstyle{remark}
  \newtheorem{remark}[subsection]{Remark}
  \newtheorem{remarks}[subsection]{Remarks}

\theoremstyle{definition}
  \newtheorem{definition}[subsection]{Definition}

\include{psfig}
\title[geometric LP ]{A geometric approach to the Littlewood-Paley
theory }
\author{Sergiu Klainerman}
\address{Department of Mathematics, Princeton University,
 Princeton NJ 08544}
\email{ seri@@math.princeton.edu}

\author{Igor Rodnianski}
\address{Department of Mathematics, Princeton University, 
Princeton NJ 08544}
\email{ irod@@math.princeton.edu}
\subjclass{35J10\newline\newline
The first author is partially supported by NSF grant 
DMS-0070696. The second author is a Clay Prize Fellow and is partially 
supported by NSF grant DMS-01007791
}
\vspace{-0.3in}
\begin{abstract}
We develop  a geometric invariant Littlewood-Paley
theory for arbitrary tensors on a compact 2 dimensional manifold. 
We show that all the important features of the classical 
LP theory survive with estimates which depend only
on very  limited regularity assumptions on the metric.
We give invariant descriptions of Sobolev and Besov spaces and 
prove some sharp  product inequalities. This theory has being developed
in connection  to the work
of the authors on the geometry of null hypersurfaces 
with a finite curvature flux condition, see \cite{KR1}, \cite{KR3}.
 We are confident however
that it can  be applied, and extended, to many
different situations.

\end{abstract} 
\maketitle
\section{introduction}
    In its simplest manifestation   Littlewood-Paley  theory
is   a systematic  method to understand 
various properties of functions $f$,
defined on ${\Bbb R}^n$,
by decomposing them in infinite  dyadic sums  $f=\sum_{k\in \Bbb Z} f_k$, 
 with frequency localized    components $f_k$, i.e. 
$\widehat{f_k}(\xi)=0$  for all values of $\xi$ outside the annulus
 $2^{k-1}\le |\xi|\le 2^{k+1}$. Such a decomposition can be  easily  achieved 
by choosing a  test function $\chi=\chi(|\xi|)$ in  Fourier 
 space,  supported in $\f12\le |\xi|\le 2$, and such that, for all $\xi\neq 0$, 
$\sum_{k\in \Bbb Z}\chi(2^{-k}\xi)=1$. Then set $\widehat{f_k}(\xi)=
\chi(2^k\xi) \hat{f}(\xi)$ or , in physical space,
$$P_k f= f_k=m_k*f$$
where $m_k(x)= 2^{nk}m(2^{k} x)$ and $m(x)$ the 
inverse Fourier transform of $\chi$. The operators $P_k$ 
are called cut-off operators or, improperly, LP projections.
We denote $P_J=\sum_{k\in J} P_k$ for all intervals $J\subset {\Bbb Z}$.

The following properties of these LP projections  are 
very easy to verify and lie  at the heart of the classical LP theory:

{\bf LP 1.}\quad  {\sl Almost Orthogonality:} \quad The operators $P_k$
are selfadjoint and verify
 $P_{k_1}P_{k_2}=0$   for all pairs of integers  such that  $|k_1-k_2|\ge 2$.
In particular,
$$\|F\|^2_{L^2}\approx\sum_k\|P_k F\|^2_{L^2}$$

{\bf LP 2.} \quad {\sl $L^p$-boundedness:} \quad For any $1\le
p\le \infty$, and any interval $J\subset \Bbb Z$,
\be{eq:pdf1'}
\|P_JF\|_{L^p}\les \|F\|_{L^p}
\end{equation}

{\bf LP 3.} \quad {\sl Finite band property:}\quad  We can write any partial derivative
$\nab P_k f$ in the form  $\nab P_k f=2^k \tilde{P}_k f$ where $\tilde{P}_k$
are the  LP-projections associated with a slightly  different test function $\tilde{\chi}$ 
and  verify the property {\bf LP2}. Thus, in particular,
for any $1\le p\le \infty$
\beaa
\|\nab P_k F\|_{L^p}&\les& 2^{k}\|  F\|_{L^p}\\
2^k\|P_k F\|_{L^p}&\les& \|\nab F\|_{L^p}
\eeaa

 {\bf LP 4.}  \quad{\sl Bernstein inequalities}. \quad For any $2\le p\le \infty$
 we have the Bernstein inequality and its dual,
$$  \|P_k F\|_{L^\infty}\les 2^{k(1-\frac 2p)} \|F\|_{L^2},\qquad 
\|P_k F\|_{L^2}\les 2^{k(1-\frac 2p)} \|F\|_{L^{p'}}$$

 The last two properties  go a long way to explain why LP theory 
is such a  useful  tool for  partial differential equations. The finite band 
property allows us to replace derivatives of the dyadic
 components $f_k$ by multiplication with $2^k$. The $L^2\to L^\infty$ Bernstein 
 inequality
is a dyadic remedy for the  failure of the 
embedding 
 of the Sobolev space $H^{\frac{n}{2}}({\Bbb R}^n)$ to $L^\infty({\Bbb R}^n)$. 
Indeed, in view
of the finite band property, the Bernstein inequality does actually imply the desired Sobolev
inequality for each LP  component $f_k$, the failure  of the Sobolev  inequality
for $f$ is due to the summation $f=\sum_k f_k$. 

Just like Fourier transform,  Littlewood-Paley theory allows one to
separate   waves of various frequencies for 
linear partial differential equations with
 constant coefficients and therefore its usefulness 
 in this context is not that surprising. 
 It took 
longer to realize that it is  helpful, in fact even more helpful,
 for the
analysis of nonlinear  equations. It turns out that multiplication
properties of various classical spaces of functions 
are   best understood by decomposing the corresponding functions
in dyadic LP components. This allows one to isolate and treat differently 
interactions of various components of the functions.
 Moreover the LP calculus allows one to manipulate a nonlinear
 PDE  to derive coupled equations for each particular frequency.  A 
first systematic application\footnote{The first manifestation 
of these type of ideas can be traced to the work of J. Nash on the isometric 
embedding problem \cite{Na}}
 of LP theory
to nonlinear PDE's was developed by Bony in the form of what is called  the 
paradifferential calculus \cite{B}. Notable applications  of LP theory include 
recent advanced in fluid dynamics, nonlinear dispersive
as well as nonlinear  wave equations (both semilinear and quasilinear),
see e.g.  \cite{Ch}, \cite{Ba-Ch}, \cite{Bour}, \cite{Tat},
 \cite{Tao}, \cite{Sm-Ta}.    

In this paper we develop an invariant LP  theory  for 
compact 2-surfaces. Our immediate goal is to apply this
theory to study the geometric properties of null
hypersurfaces, in Einstein-vacuum manifolds, with
a finite curvature flux condition, see \cite{KR1}-\cite{KR3}.  We believe however
that the theory we develop can have far wider applications.

Following a well-known procedure (see Stein \cite{S1})  we base our approach on 
heat flow,
\be{eq: heat-flow-intr}
\pr_\tau U(\tau)F -\lap U(\tau) F=0, \,\, U(0)F=F
\end{equation}
with $\lap=g^{ij}\nab_i\nab_j$ the usual Laplace-Beltrami operator  defined on the space of smooth 
tensorfields of order $m\ge 0$. 

 We then  define LP
projections $P_k$
according to the formula,
 \be{eq:LP-introd}P_k F=\int_0^\infty m_k(\tau) U(\tau) F d\tau
\end{equation}
where $m_k(\tau)=2^{2k}m(2^{2k}\tau)$ and $m(\tau)$ is a Schwartz function 
with a finite number of  vanishing moments.
 
Under some  primitive assumptions on the geometry
of our compact 2-dimensional manifold $S$     we prove  a sequence of properties for these
geometric LP  projections, similar to {\bf LP1}-- {\bf LP4}. Some of our results are
necessarily weaker\footnote{Indeed, even in Euclidean space the LP projections 
constructed by the heat flow do not possess sharp localization properties in Fourier 
space}. For example the pointwise version of the almost orthogonality property
{\bf LP1} does not hold. We can replace it however by its sufficiently
robust $L^p$  analogue.  We also  find satisfactory  analogues
for {\bf LP2}-{\bf LP4}. However we discover that the minimal geometric assumptions, 
we impose, restrict the range of $p$ in {\bf LP 3} to $p=2$ and $p\ne \infty$ in 
{\bf LP 4}.  Moreover, the $L^2\to L^\infty$ Bernstein 
inequality requires additional geometric assumptions which differ dependent 
on whether $F$ is a scalar or a tensor.

In section 2  we state our main regularity assumptions on a 2-D manifold $S$ 
and establish some basic calculus inequalities. This is the only place in the
paper where we make use of special coordinates.  Our assumption of 
{\sl weak regularity} is meant to guarantee  the existence
of such  coordinates.

Section 3  discusses the B\"ochner identities for scalar functions and 
general tensorfields. Note that the B\"ochner identity for tensorfields
has an additional term, not present for scalars, which requires 
stronger assumptions on the Gauss curvature $K$ of our manifold.

In section 4  we define the  heat flow generated by the Laplace-Beltrami 
operator $\Delta$  on tensorfields of arbitrary order. The properties 
of the heat equation derived in that section 
requires   no regularity
assumptions on $S$ beyond the fact that the metric must be Riemannian.

In section 5 we use the heat flow to develop an invariant, tensorial, 
Littlewood-Paley theory on manifolds. We prove analogues of the 
{\bf LP1}--{\bf LP4} properties of the classical LP theory.  
Once more, for most properties of our LP
projections,  we need  no regularity
assumptions on the metric, beyond the fact that
it is Riemannian.  We do  however    make  use 
of the   weak  regularity  assumption on
our  manifold $S$ in
the proof of the weak Bernstein inequality and 
its consequences. 

In sections 7   and 8  we define  fractional Sobolev and Besov spaces.

 In Section 9   we    show
how to use the geometric  LP theory developed
so far to    prove some (non sharp)    product estimates in fractional
Sobolev and Besov spaces.

In section 10   we discuss the sharp $L^2\to L^\infty$ Bernstein 
inequality. In addition to the main weak regularity assumptions 
on the 2-D manifold $S$ we have to impose conditions on its Gauss 
curvature $K$. We detect a sharp difference in the requirements imposed
on $K$ dependent on whether  we consider the scalar or the the general
tensorial  case.

In section 11  we return to the  earlier product estimates and prove
 their sharp
versions under the additional conditions needed for the sharp Bernstein 
inequality.

In section 12. we consider the mapping property of the
 covariant differentiation
$\nab$ on the Besov space $B^1_{2,1}$.   

{\bf Acknowledgments}:\,\, {\em We want to thank E. Stein for
 very helpful discussions and suggestions. }

\section{Calculus inequalities}
In this section we establish some basic calculus inequalities 
on a smooth,  compact, 2 -D maniflold $S$. We say that a coordinate
chart $U\subset M$  with coordinates $x^1, x^2$ is admissible if,
relative to  these coordinates, there exists  a constant $c>0$
such that,
\be{eq:coordchart}
c^{-1}|\xi|^2\le \ga_{ab}(p)\xi^a\xi^b\le c|\xi|^2, \qquad \mbox{uniformly for  all }
\,\, p\in U
\end{equation}
We also assume that the Christoffel symbols $\Ga^a_{bc}$ verify,
\be{eq:gammaL2}
\sum_{a,b,c}\int_U|\Ga^a_{bc}|^2 dx^1dx^2\le c^{-1}
\end{equation}
\begin{definition}
We say that a a smooth  2-d manifold $S$ is {\it weakly regular} ({\bf WR}) if 
can be covered by a finite number
of admissible coordinate charts, i.e., charts satisfying 
the conditions \eqref{eq:coordchart}, \eqref{eq:gammaL2}.
\end{definition}
 \begin{remark}
 Although we assume that our manifold $S$ is smooth  our results below
depend only on the constants in \eqref{eq:coordchart}
  and \eqref{eq:gammaL2}.
 The notion of weak regularity is introduced  to emphasize this
 fact.
 \end{remark}
Whenever 
we have inequalities of the type $A\le C\c B$,
with  C a constant which depends only on   $c$ above,   we write 
$A\les B$.

Under the {\bf WR} assumption  a it is easy
to prove the following calculus inequalities:
\begin{proposition} Let $f$ be a  real scalar function  on a 2-d weakly 
regular manifold $S$.
Then, 
\be{eq:isoperimetric}
\|f\|_{L^2}\les\|\nab f\|_{L^1}+\|f\|_{L^1}
\end{equation}
 \be{eq:LinftyL1}
\|f\|_{L^\infty}\les  \|\nab^2 f\|_{L^1}+\|f\|_{L^2}
\end{equation}
\end{proposition}
\begin{proof}:\quad Both statements  can be
reduced, by a partition of unity, to the case when the  function
 $f$   has compact support in  an admissible local chart  $U\subset S$.
Let $x^1,x^2$ be an admissible   system of   coordinates in $U $. 
 Then,
\beaa
|f(x^1, x^2)^2|&=&\big|\int_{-\infty}^{x^1}\pr_1 f(y,x^2) dy\c
 \int_{-\infty}^{x^2}\pr_2 f(x^1,y) dy\big|\\
&\les&\int_{-\infty}^{\infty}|\pr_1 f(y,x^2) dy|\c
 \int_{-\infty}^{\infty}|\pr_2 f(x^1,y)| dy
\eeaa
Hence, 
\beaa
\int_{{\Bbb R}^2}|f(x^1, x^2)|^2dx^1dx^2&\les& 
\int_{{\Bbb R}^2}|\pr_1 f(x^1, x^2)|dx^1dx^2\c\int_{{\Bbb R}^2}|\pr_2 f(x^1, x^2)|dx^1dx^2\\
&\les&\int_{{\Bbb R}^2}|\nab f(x^1, x^2)|dx^1dx^2.
\eeaa
Thus,  since in view of \eqref{eq:coordchart} $c\le \sqrt{|g|} \le c^{-1}$,
$$\big(\int_U |f(x)|^2\sqrt{|g|} dx^1dx^2\big)^\f12\les
 \big(\int_U |\nab f(x)|^2\sqrt{|g|} dx^1dx^2\big)^\f12.$$
as desired.
Simmilarly,
$$f(x^1,x^2)= \int_{-\infty}^{x^1}\int_{-\infty}^{x^2}\pr_1\pr_2 f(y^1,y^2) dy^1 dy^2.$$
Hence,
\beaa
|f(x^1,x^2)|&\le& \int_{{\Bbb R}^2}\big(|\nab^2 f(y^1,y^2)|+|\Ga||\nab f(y^1,y^2)|\big)\\
&\les&\int_S|\nab^2 f| +\big(\int_U|\Ga|^2\big)^\f12 \|\nab f\|_{L^2(S)}\\
&\les&\|\nab^2 f\|_{L^1(S)}+\|\nab f\|_{L^2(S)}
\eeaa
\end{proof}
 
As a corollary of the estimate \eqref{eq:isoperimetric}
we can derive the following Gagliardo-Nirenberg inequality:
\begin{corollary}
Given an arbitrary tensorfield  $F$ on $M$ and any $2\le p<\infty$ we have,
\be{eq:GNirenberg}
\|F\|_{L^p}\les \|\nab F\|_{L^2}^{1-\frac{2}{p}}\|F\|_{L^2}^{\frac 2 p}+\|F\|_{L^2}
\end{equation}
\end{corollary}
\begin{proof}:\quad
For any $p\ge 2$ we can write,
\beaa
\|F \|_{L^p}^{p/2}&=&\||F|^{p/2}\|_{L^2}\les\|\nab |F| ^{p/2}\|_{L^1}
+\||F|^{p/2}\|_{L^1}\\
&\les&\big(\|\nab F\|_{L^2}+\|F\|_{L^2}\big)\c
\|F\|_{L^{p-2}}^{\frac{p-2}{2}}
\eeaa
Thus, inductively, for all $p=2k$, $k=1,2,\ldots$
$$\|F\|_{L^{2k}}\les\big(\|\nab F\|_{L^2}+\|F\|_{L^2}\big)^{1-\frac{1}{k}}
\c\|F\|_{L^2}^{\frac{1}{k}}
$$
The result for general $p$ now follows by  interpolation in the scale of 
$L^p$ spaces. 
\end{proof}
As a Corollary to \eqref{eq:LinftyL1} we also  derive
\begin{corollary} For any tensorfield $F$ on $S$,
\be{eq:LinftyL2tensor}
\|F\|_{L^\infty}\les \|\nab ^2 F\|_{L^2}^\f12
\c\|F\|_{L^2}^\f12+\|F\|_{L^2}
\end{equation}
Moreover, we have a more precise estimate for any $2\le p<\infty$,
\be{eq:LinftyLp}
\|F\|_{L^\infty}\les \|\nab^2 F\|^{\frac 1p}_{L^2} 
\big (\|\nab F\|_{L^2}^{\frac {p-2}{p}} \|F\|_{L^2}^{\frac 1{p}} + 
\|F\|^{\frac {p-1}p}_{L^2}\big ) 
+ \|\nab F\|_{L^2}.
\end{equation}

\end{corollary}
\begin{proof}:\quad  We apply the estimate \eqref{eq:LinftyL1} to the scalar $|F|^2$
as follows,
\beaa
\|F\|_{L^\infty}^2&\les& \|\nab^2 |F|^2\|_{L^1}+\||F|^2\|_{L^2}\\
&\les&\|\nab^2 F\|_{L^2}\|F\|_{L^2}+\|\nab F\|_{L^2}^2+ \|F\|_{L^4}^2\\
\eeaa
In view of \eqref{eq:GNirenberg}, 
\beaa
\|F\|_{L^4}^2\les \|\nab F\|_{L^2}\|F\|_{L^2}+\|F\|_{L^2}^2
\eeaa
Hence,
\beaa
\|F\|_{L^\infty}^2\les\|\nab^2 F\|_{L^2}\|F\|_{L^2}+\|\nab F\|_{L^2}^2+
\|\nab F\|_{L^2}\|F\|_{L^2}+\|F\|_{L^2}^2
\eeaa
The desired estimate now follows by Cauchy-Schwartz.
To prove the estimate \eqref{eq:LinftyLp} we observe that applying 
\eqref{eq:LinftyL1} to $|F|^p$ we obtain
$$
\|F\|_{L^\infty}\les \|\nab^2 F\|^{\frac 1p}_{L^2} \|F\|_{L^{2(p-1)}}^{\frac {p-1}p} + \|\nab F\|_{L^2}
$$

By the Galgiardo-Nirenberg inequality \eqref{eq:GNirenberg} we have that 
$$
\|F\|_{L^{2(p-1)}}\les \|\nab F\|_{L^2}^{\frac {p-2}{p-1}} \|F\|_{L^2}^{\frac 1{p-1}}
+ \|F\|_{L^2}
$$
Thus, finally
$$
\|F\|_{L^\infty}\les \|\nab^2 F\|^{\frac 1p}_{L^2} 
\big (\|\nab F\|_{L^2}^{\frac {p-2}{p}} \|F\|_{L^2}^{\frac 1{p}} +
\|F\|^{\frac {p-1}p}_{L^2}\big )
+ \|\nab F\|_{L^2}
$$
as desired.
\end{proof}

\section{B\"ochner identity}
In this section  we recall  the B\"ochner identity on a 2-D manifold.
This allows us  to 
control  the $L^{2}$ norm of the second derivatives of a tensorfield 
in terms of the $L^{2}$ norm of the laplacian and 
geometric quantities associated with a given 2-surface.
\begin{proposition}
Let  $K$ denote the  Gauss curvature of our 2-D riemannian manifold $M$.
Then

\noindent 
{\bf i})\quad For a scalar function $f$
\begin{equation}
\label{eq:scalBoch}
\int_{S} |\nab^{2} f|^{2} = \int_{S} |\lap f|^{2} - 
\int_{S} K |\nab f|^{2}
\end{equation}

\noindent
{\bf ii)}\quad For a vectorfield $F_{a}$
\begin{equation}
\label{eq:vectorBoch}
\int_{S} |\nab^{2} F|^{2} = \int_{S} |\lap F|^{2} -
\int_S K (2\,|\nab F|^{2}-|\div F|^{2}-|\curl F|^2) + \int_{S} K^{2} |F|^{2}
\end{equation}
where $\div F=\ga^{ab}\nab_b F_a$, $\curl F=\div( \dual F)=\in_{ab}\nab_a F_b$
\label{prop:Bochner}
\end{proposition}
\begin{proof}:\quad
Recall that on a 2-surface the Riemann tensor
\be{eq:curv-2d}
R_{abcd} = (\ga_{ac}\ga_{bd}- \ga_{ad}\ga_{bc})K,\qquad  R_{ab}=\ga_{ab} K,
\end{equation}
\noindent

 To prove {\bf i)} observe that, relative to an arbitrary orthonormal frame $(e_a)_{a=1,2}$,
\beaa
\nab_a (\lap f) &=&\nab_a (\nab_c\nab_c f)=\nab_c\nab_a\nab_c f+[\nab_a,\nab_c]\nab_c f\\
&=&\nab_c\nab_c \nab_a f+R_{cdac}\nab_d f\\
&=&\lap (\nab_a f)-R_{ad}\nab_d f
\eeaa 
Thus,
\beaa
\int_S|\lap f|^2&=&-\int_S \nab_a (\lap f)\c\nab_a f=\int_S \lap \nab_a f\c \nab_a f
-R_{ab}\nab_a f\nab_b f\\
&=&\int_S |\nab^2 f|^2- 
\int_{S} K |\nab f|^{2}
\eeaa
as desired.

Similarly for a vector $F_i$,
\beaa
\nab_a (\lap F_i) &=&\nab_a (\nab_c\nab_c F_i)=\nab_c\nab_a\nab_c F_i+[\nab_a,\nab_c]\nab_c F_i\\
&=&\nab_c\nab_c \nab_a f_i+\nab_c([\nab_a,\nab_c] F_i)  +  
 R_{cdac}\nab_d F_i+R_{idac}\nab_c F_d\\
&=&\lap (\nab_a F_i)+\nab_c\big(  R_{idac}F_d   \big) +  
 R_{cdac}\nab_d F_i+R_{idac}\nab_c F_d\\
&=&\lap (\nab_a F_i)+\nab_c\big(  R_{idac}F_d  \big)-R_{da}\nab_d F_i+R_{idac}\nab _c F_d
\eeaa 
Hence,
\beaa
-\int_S|\lap F|^2&=&\int_S\nab_a(\lap F_a)\\
&=&-\int_S|\nab^2F|^2-\int_S R_{idac}F_d \nab_c\nab_a F_i\\
&-&\int_S R_{da}\nab_d
F_i\nab_a F_i +\int_S R_{idac}\nab _c F_d\nab_a F_i
\eeaa
Now observe that,
\beaa
\int_S R_{idac}F_d \nab_c\nab_a F_i=\f12\int_S R_{idac}F_d( \nab_c\nab_a F_i-\nab_a\nab_c F_i)
=\f12\int_S R_{idac} R_{imca} F_d F_m
\eeaa
Therefore,
\beaa
\int_S|\lap F|^2&=&\int_S|\nab^2F|^2+\f12\int_S R_{diac} R_{miac}   F_d  F_m
+\int_S R_{da}\nab_d F_i\nab_a F_i\\
&-&\int_S R_{idac}\nab _c F_d\nab_a F_i
\eeaa
Using the formulas \eqref{eq:curv-2d}  and observing that $ \nab_b F_a-\nab _a F_b=\in_{ba}
\curl F$ we find,
\beaa
R_{diac} R_{miac} F_d F_m& =&2K^2 \de_{dm}F_d F_m=2K^2|F|^2\\
R_{da}\nab_d F_i\nab_a F_i&=& K |\nab F|^2\\
R_{idac}\nab _c F_d\nab_a F_i&=&K\big(|\div F|^2-\nab_a F_b\nab_b F_a\big)\\
&=&
K\big(|\div F|^2-|\nab F|^2-\nab_a F_b(\nab_b F_a-\nab_a F_b)\big)\\
&=&K\big(|\div F|^2+|\curl F|^2-|\nab F|^2\big)
\eeaa
Therefore,
\beaa
\int_S|\lap F|^2=\int_S|\nab^2F|^2-\int_S
|K|^2 |F|^2+\int_SK\big((2\,|\nab F|^2-(|\div F|^2+|\curl F|^2)\big)
\eeaa
as desired.
\end{proof}
\begin{corollary}[B\"ochner inequality]
For any tensorfield $F$ and an arbitrary $2\le p<\infty$
\bea
\|\nab^2 F\|_{L^2} &\les &\|\Delta F\|_{L^2} + 
( \|K\|_{L^2} + \|K\|_{L^2}^\f12)\|\nab F\|_{L^2}\label{eq:Bochner-ineq}\\
&+&
 \|K\|_{L^2}^{\frac p{p-1}}
\big ( \|\nab F\|_{L^2}^{\frac {p-2}{p-1}} \|F\|_{L^2}^{\frac 1{p-1}} 
 + \|F\|_{L^2}\big )
\eea
\label{cor:Bochner-Ineq}
\end{corollary}
\begin{proof}:\quad 
The B\"ochner identity \eqref{eq:vectorBoch} implies that
\be{eq:Bochner1}
\|\nab^2 F\|_{L^2} \les \|\Delta F\|_{L^2} + \|K\|_{L^2}^\f12
\|\nab F\|_{L^4} + \|K\|_{L^2} \|F\|_{L^\infty}
\end{equation}
Using the Gagliardo-Nirenberg inequality \eqref{eq:GNirenberg}
and the estimate \eqref{eq:LinftyLp} we infer that for any $2\le p<\infty$
\begin{align*}
&\|\nab F\|_{L^4} \les \|\nab^2 F\|_{L^2}^\f12 \|\nab F\|_{L^2}^\f12
+ \|\nab F||_{L^2},
\\ &
\|F\|_{L^\infty}\les \|\nab^2 F\|^{\frac 1p}_{L^2}
\big ( \|\nab F\|_{L^2}^{\frac {p-2}{p}} \|F\|_{L^2}^{\frac 1{p}} 
+ \|F\|_{L^2}^{\frac {p-1}p}\big )
+ \|\nab F\|_{L^2}
\end{align*}
Substituting this into \eqref{eq:Bochner1} we obtain
\beaa
\|\nab^2 F\|_{L^2} &\les & \|\Delta F\|_{L^2} +
 \|K\|_{L^2}^\f12\Big ( \|\nab^2 F\|_{L^2}^\f12 \|\nab F\|_{L^2}^\f12 
 + \|\nab F\|_{L^2}\Big ) \\ &+&
 \|K\|_{L^2} 
\Big (\|\nab^2 F\|^{\frac 1p}_{L^2} \big (\|\nab F\|_{L^2}^{\frac {p-2}{p}} \|F\|_{L^2}^{\frac 1{p}} + \|F\|_{L^2}^{\frac {p-1}p}\big ) + \|\nab F\|_{L^2}\Big )
\eeaa
This, in turn, implies that 
$$
\|\nab^2 F\|_{L^2} \les \|\Delta F\|_{L^2} + 
( \|K\|_{L^2} + \|K\|_{L^2}^\f12)\|\nab F\|_{L^2}+
 \|K\|_{L^2}^{\frac p{p-1}}
\big ( \|\nab F\|_{L^2}^{\frac {p-2}{p-1}} \|F\|_{L^2}^{\frac 1{p-1}} 
 + \|F\|_{L^2}\big )
$$
 as desired.
 \end{proof}

\section{Heat equation on $S$}
In this section we study the properties of the heat  equation  for arbitrary
tensorfields $F$ on $S$.
$$\pr_\tau U(\tau)F -\lap U(\tau) F=0, \,\, U(0)F=F,$$
with $\lap=\lap_\ga$ the usual Laplace-Beltrami
operator on $S$. Observe that  the operators $U(\tau)$ are
selfadjoint\footnote{Indeed observe that
$\lap$  is selfadjoint and formally $U(\tau) f =\sum_n \frac{1}{n!} t^n\lap^n$.} 
and form a semigroup for $\tau>0$. In other words  for all, real valued, smooth 
tensorfields $F, G$,
\be{eq:seladj-semigroup}
\int_SU(\tau)F\c G =\int_SF\c U(\tau) G, \qquad U(\tau_1)U(\tau_2)=U(\tau_1+\tau_2)
\end{equation}
We shall prove the following $L^2$ estimates for the operator
$U(\tau)$.
\begin{proposition}
We have the following estimates for the operator $U(\tau)$:
\begin{align}
&\|U(\tau) F\|_{L^2(S)}\le \|F\|_{L^2(S)}\label{eq:l2heat1}\\
&\|\nab U(\tau) F\|_{L^2(S)}\le \|\nab F\|_{L^2(S)}\label{eq:l2heatnab}\\
&\|\nab U(\tau) F\|_{L^2(S)}\le \frac{\sqrt
2}{2}\tau^{-\f12}\|F\|_{L^2(S)}\label{eq:l2heat2}
\\ &\|\lap U(\tau) F\|_{L^2(S)}\le \frac{\sqrt{2}}{2}\tau^{-1}
\|F\|_{L^2(S)}
\label{eq:l2heat3}
\end{align}
\label{le:L2heat}
We also have,
\be{eq:l2heat4}
\| U(\tau)\nab F\|_{L^2(S)}\le \frac{\sqrt{2}}{2}\tau^{-\f12}\|F\|_{L^2(S)}
\end{equation}
\end{proposition}
\begin{proof}:\,\,\,
To prove \eqref{eq:l2heat1} we multiply the equation
$$
\pr_\tau U(\tau)F  - \lap U(\tau) F =0
$$
by $U(\tau) F$ and integrate over $S$.
$$
\frac 12\frac {d}{d\tau} \|U(\tau) F\|_{L^2(S)}^2+ \|\nab U(\tau)F\|_{L^2(S)}^2=0
$$
Therefore,
\begin{equation}
\label{eq:ident1}
\frac 12\|U(\tau) F\|_{L^2(S)}^2+ \int_0^\tau\|\nab U(\tau')F\|_{L^2(S)}^2 d\tau'=
\frac 12\|F\|^2_{L^2(S)}
\end{equation}
and \eqref{eq:l2heat1} follows.
On the other hand, multiplying the equation by $\tau \lap U(\tau) F$, we similarly 
obtain the identity
$$
\frac 12\frac {d}{d\tau} \tau
\|\nab U(\tau) F\|_{L^2(S)}^2+ \tau \|\lap U(\tau)F\|_{L^2(S)}^2= 
\frac 12\|\nab U(\tau)F\|_{L^2(S)}^2
$$
 Integrating this in $\tau$, with the help of \eqref{eq:ident1}, 
\begin{equation}
\label{eq:ident2}
\frac \tau 2\|\nab U(\tau) F\|_{L^2(S)}^2+ 
\int_0^\tau \tau'\|\lap U(\tau')F\|_{L^2(S)}^2
d\tau'\le \f12 \int_0^\tau\|\nab U(\tau)F\|_{L^2(S)}\le \frac 1 4  \|F\|^2_{L^2(S)}
\end{equation}
which implies \eqref{eq:l2heat2}. Proceeding in exactly the same way 
with the multiplier $\tau \lap U(\tau) F$ replaced by $ \lap U(\tau) F$
yields \eqref{eq:l2heatnab}.
Furthermore, multiplying the equation by $\tau^2 \lap^2 U(\tau) f$, we have
$$
\frac 12\frac {d}{d\tau} \tau^2
\|\lap U(\tau) F\|_{L^2(S)}^2+ \tau^2 \|\nab\lap U(\tau)F\|_{L^2(S)}^2= 
\tau \|\lap U(\tau)F\|_{L^2(S)}^2
$$
Integrating in $\tau$ and using \eqref{eq:ident2}, we obtain
$$
\frac {\tau^2} 2\|\lap U(\tau) F\|_{L^2(S)}^2+ 
\int_0^\tau (\tau')^2\|\nab\lap U(\tau')F\|_{L^2(S)}^2
d\tau'=\int_0^\tau \tau' \|\lap U(\tau')F\|_{L^2(S)}^2 d\tau' \le \frac 1 4\|F\|^2_{L^2(S)}
$$
This immediately yields \eqref{eq:l2heat3}.

To prove \eqref{eq:l2heat4} we observe that
$$\|U(\tau)\nab F\|_{L^2}^2=<U(\tau)\nab F\,,\, U(\tau)\nab F>=< \div
U(\tau)U(\tau)\nab F\,,\,  F>$$ Therefore,
\beaa
\|U(\tau)\nab F\|_{L^2}^2&\le&
 \|\nab (U(\tau)U(\tau)\nab F)\|_{L^2}\|F\|_{L^2}\\
&\le&\frac{\sqrt{2}}{2}\tau^{-\f12}\|U(\tau)\nab F\|_{L^2}\|F\|_{L^2}
\eeaa
whence $\|U(\tau)\nab F\|_{L^2}\les \frac{\sqrt{2}}{2}\tau^{-\f12}\|F\|_{L^2}$ as desired.
\end{proof}
In the next proposition
we establish a simple $L^p$ estimate for $U(\tau)$.
\begin{proposition} For every $2\le p\le\infty$
  $$ \|U(\tau)F\|_{L^p}\le \|F\|_{L^p}$$
\end{proposition}
 
\begin{proof}:\quad  We shall first prove
the Lemma for scalar functions $f$.
 We multiply the equation $\pr_\tau U(\tau)f- \lap
U(\tau)f=0$ by $\big(U(\tau)f\big)^{2p-1}$ and integrate by parts.
We get,
$$\frac{1}{2p}\frac{d}{d\tau} \|U(\tau)F\|_{L^{2p}}^{2p}+
(2p-1)\int |\nab U(\tau) f|^2|U(\tau)f|^{2p-2}=0$$
Therefore,
$$\|U(\tau) F\|_{L^{2p}}\le \|F\|_{L^{2p}}$$
The case when  $F$ is a tensorfield can be treated in the same manner 
with  multiplier   $\big(|U(\tau)F|^2\big)^{p-1}U(\tau) F$.
\end{proof}
\section{ Invariant Littlewood-Paley theory}
In this section we shall
use the heat flow discussed in the previous section
to develop  an invariant, fully tensorial, Littlewood-Paley theory on manifolds.
Though we restrict ourselves here to two dimensional compact  manifolds
it is clear that our theory can be extended to arbitrary dimensions and noncompact
manifolds.
\begin{definition}
Consider  the class $\cal M$ of smooth functions $m$ on $[0,\infty)$,
vanishing sufficiently fast at $\infty$,
verifying the  vanishing  moments property:
\be{eq:moments}
\int_0^\infty \tau^{k_1}\pr_\tau^{k_2} m(\tau) d\tau=0, \,\,\,\,
|k_1|+|k_2|\le N 
\end{equation}

 We set,
  $m_k(\tau)=2^{2k}m(2^{2k}\tau)$ 
and  define the geometric Littlewood -Paley (LP) 
projections $P_k$, associated to the LP- representative 
function $m\in \MM$, for   arbitrary tensorfields  $F$ on $S$
to be 
\be{eq:LP}P_k F=\int_0^\infty m_k(\tau) U(\tau) F d\tau
\end{equation}
Given an interval $I\subset \Bbb Z$ we define $$P_I=\sum_{k\in I} P_k F.$$
In particular we shall use the notation $P_{<k}, P_{\le k}, P_{>k}, P_{\ge k}$.
\end{definition}
Observe that $P_k$ are selfadjoint\footnote{This follows easily 
in view of the selfadjoint properties
of $\Delta$ and $U(\tau)$.}, i.e., $P_k=P_k^*$, in the sense,
$$<P_kF, G>=<F,P_k G>,$$
where, for any given $m$-tensors $F,G$ 
$$<F,G>=\int_S\ga^{i_1j_1}\ldots\ga^{i_mj_m}
F_{i_1\ldots i_m}G_{j_1\ldots j_m}d\mbox{vol}_\ga    $$ 
denotes the usual $L^2$ scalar product.

 Consider two LP
projections associated to $a, b$
\beaa
P_aP_b F&=&
\int_0^\infty\int_0^\infty d\tau_1 d\tau_2 a(\tau_1)b(\tau_2) U(\tau_1+\tau_2)F\\
&=&\int_0^\infty d\tau_1 \int_{\tau_1}^\infty d\tau \,\,a(\tau_1)b(\tau-\tau_1)
U(\tau)F\\
&=&\int_0^\infty d\tau U(\tau)f\int_{0}^\tau d\tau_1 \,\,a(\tau_1)b(\tau-\tau_1)\\
&=&\int_0^\infty d\tau \,a\star b(\tau)\,   U(\tau)F
\eeaa
where 
\be{eq:funnyconv}
a\star b(\tau)=\int_{0}^\tau d\tau_1 \,\,a(\tau_1)b(\tau-\tau_1)
\end{equation}
\begin{lemma}
If $a,b\in {\cal M}$  so does $a\star b$. 
Also, $(a\star b)_k=a_k\star b_k$. In particular
if we denote by $^{(a)}P_k$ and $^{(b)}P_k$ the LP projections
associated to  $a,b$ then,
$$^{\,(a)}P_k \c ^{\,(b)}P_k=^{\,(a\star b)}P_k$$
\label{le:convolution}
\end{lemma}
\begin{proof}:\quad We only need to show that
$\int (a\star b)(\tau) d\tau=0$. Then, 
we can easily check that $\tau\c(a\star b)(\tau)$ and $\frac{d}{d\tau}(a\star b)$
also verify the same property as well as any combination of these. Clearly
$\int_0^\infty a\star b\,
d\tau=\int_0^\infty a(\tau_1)d
\tau_1
\c
\int_0^\infty b(\tau_2) d\tau_2 =0$.
\end{proof}
Motivated by this Lemma
we define:
\begin{definition} Given a positive integer $\ell$ we 
 define the class $\Ml\subset\MM$ of LP-
representatives 
 to consist of functions of the form 
$$\bar{m}=m\star m\star\ldots\star m=(m\star)^\ell,$$
for some $m\in \MM$.\label{def:Ml}
\end{definition}
\begin{lemma} For any $\ell\ge 1$  there exists an element $\bar m\in \Ml$ 
such that the LP-projections associated to $\bar m$ verify:
\be{eq:partition}
\sum_kP_k=I
\end{equation}
\end{lemma}
\begin{proof}
!!!
\end{proof}
\begin{theorem} The LP-projections $P_k$ associated
to an arbitrary $m\in {\cal M}$ verify the following
 properties:

i)\quad {\sl $L^p$-boundedness} \quad For any $1\le
p\le \infty$, and any interval $I\subset \Bbb Z$,
\be{eq:pdf1}
\|P_IF\|_{L^p}\les \|F\|_{L^p}
\end{equation}

ii) \quad {\sl $L^p$- Almost Orthogonality}\quad  Consider two families
of LP-projections $P_k, \tilde P_k$ associated to $m$ and  respectively 
$\tilde m$, both in ${\cal M}$. For any  
$1\le p\le
\infty$:
\be{eq:pdf2}
\|P_k\tilde P_{k'}F\|_{L^p}\les 2^{-2|k-k'|} \|F\|_{L^p}
\end{equation}

iii) \quad  {\sl Bessel inequality} 
$$\sum_k\|P_k F\|_{L^2}^2\les \|F\|_{L^2}^2$$
iv)\quad {\sl Reproducing Property} \quad  Given any integer $\ell\ge 2$ 
 and  $\bar m\in \Ml$ there exists $m\in \MM$ such that 
 such that $\bar m=  m\star m$. Thus,
$$^{(\bar m)}P_k =^{(m)}P_k\c ^{(m)}P_k.$$
Whenever there is no danger of confusion we shall simply write $P_k=P_k\c P_k$.

v)\quad {\sl Finite band property}\quad For any $1\le p\le \infty$.
\beaa
\|\lap P_k F\|_{L^p}&\les& 2^{2k} \|F\|_{L^p}\\
\|P_kF\|_{L^p} &\les& 2^{-2k} \|\lap F \|_{L^p}
\eeaa
Moreover given $m\in \MM$ we can find $\bar{m}\in \MM$
such that $\lap P_k =2^{2k}\bar{P}_k$ with $\bar{P}_k$
the LP projections associated to $\bar{m}$.

In addition, the $L^2$ estimates
\beaa
\|\nab P_k F\|_{L^2}&\les& 2^{k} \|F\|_{L^2}\\
\|P_kF\|_{L^2} &\les& 2^{-k} \|\nab F  \|_{L^2}
\eeaa
hold together with the dual estimate
$$\| P_k \nab F\|_{L^2}\les 2^k \|F\|_{L^2}$$

vi) \quad{\sl Weak Bernstein inequality}\quad For any $2\le p<\infty$
\begin{align*}
&\|P_k F\|_{L^p}\les (2^{(1-\frac 2p)k}+1) \|F\|_{L^2},\\
&\|P_{<0} F\|_{L^p}\les \|F\|_{L^2}
\end{align*}
together with the dual estimates 
\begin{align*}
&\|P_k F\|_{L^2}\les (2^{(1-\frac 2p)k}+1) \|F\|_{L^{p'}},\\
&\|P_{<0} F\|_{L^2}\les \|F\|_{L^{p'}}
\end{align*}
vii)\quad{\sl Commutator Estimate} \quad Given two tensorfields
$F, G $ and  $F\c G$ any contraction of the tensor product  $F\otimes G$ 
we  have the following  estimate for  the commutator
  $[P_k\,, \,F]\c G= P_k(F\c G)-F\c P_kG $
$$\|\,\,[P_k\,,\, F]\c G\|_{L^2}\les 2^{-k} \|\nab F\|_{L^\infty}\|G\|_{L^2}.$$
We also have the estimate of the form 
$$\|\,\,[P_k\,,\, F]\c G\|_{L^2}\les \bigg (
2^{-2k}\|\lap F\|_{L^{\infty}}+ 2^{-k} \|\nab F\|_{L^\infty}\bigg )
\|G\|_{L^2}.$$
\label{thm:LP}
\end{theorem}
\begin{proof}:\quad 

i) \quad The $L^p$ boundedness of $P_k$ follows from the $L^p$ mapping 
properties of the heat flow $U(\tau)$.

ii)\quad  Assume that $k_2\ge k_1$. By definition and in view of the semigroup property of
$U(\tau)$ we write,
$$P_{k_1}\tilde P_{k_2}f=\int_0^\infty\int_0^\infty U(\tau_1+\tau_2)f\,\,
m_{k_1}(\tau_1)\tilde m_{k_2}(\tau_2)d\tau_1 d\tau_2$$
Writing $U(\tau_1+\tau_2)=U(\tau_1)+\int_0^1 \frac{d}{ds}U(\tau_1+s\tau_2) ds$
and then using the vanishing of $\int_0^\infty \tilde m_{k_2} $
we infer that,
\beaa
P_{k_1}\tilde P_{k_2}f&=&\int_0^\infty\int_0^\infty \frac{d}{d\tau_1}
\int_0^1 U(\tau_1+s\tau_2)f\,\,
m_{k_1}(\tau_1)\,\tau_2\tilde m_{k_2}(\tau_2)d\tau_1 d\tau_2\\
&=&-\int_0^\infty\int_0^\infty 
\int_0^1 U(\tau_1+s\tau_2)f\,\,
\frac{d}{d\tau_1}m_{k_1}(\tau_1)\,\tau_2\tilde m_{k_2}(\tau_2)d\tau_1 d\tau_2\\
&-&m_{k_1}(0)\int_0^\infty d\tau_2 \,\tau_2\tilde m_{k_2}(\tau_2)\int_0^1 U(s\tau_2) f ds
\eeaa
Now setting $ \tilde n(\tau)=\tau \tilde m(\tau) $,  and $n(\tau)=m'(\tau)$
we infer that,
\beaa
P_{k_1}\tilde P_{k_2}f&=&-2^{2(k_1-k_2)}\int_0^\infty\int_0^\infty 
\int_0^1 U(\tau_1+s\tau_2)f\,\,
n_{k_1}(\tau_1)\,\tilde n_{k_2}(\tau_2)d\tau_1 d\tau_2\\
&-&2^{2(k_1-k_2)}m(0)\int_0^\infty d\tau_2 \tilde n_{k_2}(\tau_2)\int_0^1 U(s\tau_2) f
ds
\eeaa
Therefore,
using the $L^p$ mapping properties of $U$,
 \beaa
\|P_{k_1}\tilde P_{k_2}F\|_{L^p}&=&2^{-2|k_1-k_2|}\|F\|_{L^p}\int_0^\infty\int_0^\infty 
\,\,
|n_{k_1}(\tau_1)|\,|\tilde n_{k_2}(\tau_2)|d\tau_1 d\tau_2\\
&+&2^{-2|k_1-k_2|}|m(0)|\|F\|_{L^p}\int_0^\infty|\tilde n_{k_2}(\tau_2)|d\tau_2\\
&\les&2^{-2|k_1-k_2|}\|F\|_{L^p}
\eeaa
\begin{remark}
One can   give a slicker proof of the almost orthogonality properties of 
LP projections  by using the algebraic formula $2^{2k} P_k f=\lap \bar{P}_k f$,
see \eqref{eq:Pk=2-2klapbarPk} below. Moreover,   if  sufficiently many moments 
of $m$ are zero, s.t $\tau^{2j}m, \tau^{2j}\tilde{m}$ are good symbols,
then in fact,
\be{eq:stronger-orthog}
\|P_{k_1} \tilde{P}_{k_2} F\|_{L^2}\les 2^{-2j|k_1-k_2|}\|F\|_{L^2}
\end{equation}
\label{rem:inside-thm}
\end{remark}

iii)\quad  To prove the Bessel type inequality  
we write,
$$\sum_k\|P_k F\|_{L^2}^2=\sum_k<P_kP_k f, f>\le
 \|(\sum_k P_k^2) F\|_{L^2}\|F\|_{L^2}$$
To show that the operator $P=\sum_k P_k^2$ is bounded on $L^2$ 
we appeal to the Cotlar-Stein Lemma, see \cite{S2}. Observe first that, in view
of Lemma \ref{le:convolution}, $P_k^2= ^{\,\,(m\star m)}P_k$.  
Since $m\star m\in {\cal M}$
we can, without loss of generality, simply write $P_k^2=P_k$.  
The conditions of applicability of the Cotlar-Stein Lemma\footnote{Notice
that we are in the special case of commuting selfadjoint 
operators.} 
 are   satisfied in view of the almost orthogonality established
in  part ii) as well as $P_k=P_k^*$.

iv)\quad The proof is immediate in view of the definition \ref{def:Ml}.

v)\quad According to the definition of $P_kf$ we 
have 
\beaa
\lap P_k f&=&\int_0^\infty m_k(\tau) \lap U(\tau)f
=\int_0^\infty m_k(\tau) \frac{d}{d\tau} U(\tau)f\\
&=&- m_k(0)U(0) f- \int_0^\infty \frac{d}{d\tau}m_k(\tau)  U(\tau)f\\
&=&- 2^{2k}\bigg(m(0) f+ \int_0^\infty (m')_k(\tau)  U(\tau)f\bigg)\\
\eeaa
In view of the $L^p$ properties of $U(\tau)f$ and the obvious  bound
$ \int_0^\infty |(m')_k(\tau)| d\tau \les 1$,
\be{eq:lapestimate}
 \|\lap P_k F\|_{L^p}\les 2^{2k} \|F\|_{L^p}
\end{equation}
To prove the second estimate we  introduce 
$\bar{m}(\tau)= -\int_\tau^\infty  m(\tau)$
such that $\frac{d}{d\tau}\bar{m}=m(\tau)$ 
and $\int_0^\infty |\bar{m}(\tau)| d\tau <\infty$. Observe also 
that $\bar{m}(0)=0$.
Set also, $$\bar{m}_k(\tau)= 2^{2k} \bar{m}(2^{2k}\tau)$$
\bea
2^{2k} P_kf&=&\int_0^\infty  2^{2k} m_k(\tau)  U(\tau)f
=\int_0^\infty \frac{d}{d\tau}\bar{m}_k(\tau)  U(\tau)f\nn\\
&=& -\int_0^\infty \bar{m}_k(\tau) \frac{d}{d\tau} U(\tau)f
=- \int_0^\infty \bar{m}_k(\tau) \lap U(\tau)f\nn\\
&=&- \int_0^\infty \bar{m}_k(\tau)  U(\tau)\lap f\label{eq:formula} 
\eea
Therefore, using the estimate $\| U(\tau)\lap F\|_{L^p}\les \|\lap F\|_{L^p}$,
we infer that,
\be{eq:22klp}
2^{2k} \|P_kF\|_{L^p}\les \|\lap F\|_{L^p} \int_0^\infty |\bar{m}_k(\tau) |d\tau
\les \|\lap F\|_{L^p}
\end{equation}
Observe also that, according to \eqref{eq:formula} we have
\be{eq:Pk=2-2klapbarPk}
 2^{2k}P_k F =\lap \bar{P}_k F
\end{equation}
where ${\bar P}_k$ is defined by  the  symbol $\bar{m}(\tau) =
-\int_\tau^\infty m(\tau')\in \MM.
$

To prove the $L^2$ estimates involving one derivative we observe that
\beaa
\|\nab P_k F\|_{L^2}^2 &=& <\nab P_k F, \nab P_k F> = -<\lap P_k F, P_k F>\\ &\le&
\|\lap P_k F\|_{L^2} \|P_k F\|_{L^2} \les 2^{2k} \|F\|_{L^2}^2
\eeaa
On the other hand, using \eqref{eq:formula}
\beaa
&&2^{2k}\|P_k F\|_{L^2}^2 = 2^{2k} < P_k F,  P_k F> = -\int_0^\infty \bar m_k(\tau) 
<\lap U(\tau)F, P_k F>\\ &=& \int_0^\infty \bar m_k(\tau) 
<\nab U(\tau)F,\nab P_k F>\,\,\le  \int_0^\infty |\bar m_k(\tau)|\,\c
\|\nab U(\tau)F\|_{L^2}\|\nab P_k F\|_{L^2}
\les \|\nab F\|_{L^2}^2,
\eeaa
where we used the inequality \eqref{eq:l2heatnab},
$\|\nab U(\tau)F\|_{L^2}\le \|\nab F\|_{L^2}$ together with the 
bound $\|\nab P_k F\|_{L^2}\les \|\nab F\|_{L^2}$, which follows from it.

vi) \quad 
 The proof of the $L^p$ Bernstein inequality is an easy consequence of the 
Gagliardo-Nirenberg  inequality \eqref{eq:GNirenberg}:
\be{eq:GN}
\|P_k F\|_{L^p}\les \|\nab P_k F\|_{L^2}^{1-\frac{2}{p}}\|P_k F\|_{L^2}^{\frac{2}{p}} +
\|P_k F\|_{L^2}
\end{equation}
for $2\le p<\infty$ and the finite band property.

vii) By definition 
$$[P_k \,,\, F]G=\int_0^\infty \bigg( U(\tau)(F\c G)-F\c 
U(\tau) G\bigg)m_k(\tau) d\tau$$
Let $w=U(\tau)(F\c G)-F\c U(\tau)G$. Clearly,
\beaa
 \pr_\tau w-\lap w&=&\nab\big( \nab F\c U(\tau)G \big) +\nab F\c \nab U(\tau) G\\
&=& \lap F\c U(\tau)G  +2\nab F\c \nab U(\tau) G
\eeaa
Consequently, since $w(0)=0$,
\beaa w&=&w_1+w_2\\
 w_1(\tau)&=&\int_0^\tau U(\tau-\tau')\big( \lap F\c U(\tau')G\big)  d\tau'\\
 w_2(\tau)&=&\int_0^\tau U(\tau-\tau')\big(\nab F\c \nab U(\tau') G\big)d\tau'
\eeaa
and,
\beaa
\|w_1(\tau)\|_{L^2}&\les&\int_0^\tau \| U(\tau-\tau')\big( \lap F\c U(\tau')G\big)
\|_{L^2}d\tau'
\les\int_0^\tau \|( \lap F\c U(\tau')g
\|_{L^2}d\tau'\\
&\les&\| \lap F\|_{L^\infty} \int_0^\tau \| U(\tau')G\|_{L^2}\les 
\tau\| \lap F\|_{L^\infty} 
\c\|G\|_{L^2}
\eeaa
\beaa
\|w_2(\tau)\|_{L^2}&\les&\int_0^\tau \| U(\tau-\tau')\big( \nab F\c \nab U(\tau')G\big)
\|_{L^2}d\tau'
\les\int_0^\tau \|( \nab F\c\nab  U(\tau')G
\|_{L^2}d\tau'\\
&\les& \| \nab F\|_{L^\infty} \int_0^\tau \|\nab U(\tau')G\|_{L^2}\les 
\|\nab  F\|_{L^\infty}\c\|G\|_{L^2}\int_0^\tau\tau^{-\f12} \\
&\les&\tau^\f12 \|\nab  F\|_{L^\infty}\c\|G\|_{L^2}
\eeaa
Therefore,
\beaa
\|[P_k \,,\, F]G\|&\les &\int_0^\infty \|w(\tau)\|_{L^2} |m_k(\tau)| d\tau\\
&\les&\bigg( 2^{-2k}\| \lap F\|_{L^\infty} +2^{-k}\|\nab  F\|_{L^\infty}\bigg)
\c\|G\|_{L^2}
\eeaa
\begin{remark}
To get the inequality
$$\|[P_k \,,\, F]G\|\les2^{-k}\|\nab  F\|_{L^\infty}\|G\|_{L^2}
$$
we need the $L^2$  estimate
$\|U(\tau)\nab F\|_{L^2}\les \tau^{-\f12} \|F\|_{L^2}$
established in \eqref{eq:l2heat4}.
We rewrite 
\begin{align*}
& w_1(\tau)= w_{11}(\tau)-w_{12}(\tau),\\
&w_{11}(\tau)= 
\int_0^\tau U(\tau-\tau')\nab \big( \nab F\c U(\tau')G\big)  
d\tau',\\
&w_{12}(\tau)= 
\int_0^\tau U(\tau-\tau')\big( \nab F\c \nab U(\tau')G\big)  
d\tau'
\end{align*}
The term $w_{12}$ is exactly the same as $w_{2}(\tau)$ and gives 
rise to the desired estimate. To estimate $w_{11}$ we use 
\eqref{eq:l2heat4} and write
$$
\|w_{11}(\tau)\|_{L^{2}}\les \int_{0}^{\tau} {\tau'}^{-\f12} \|\nab F\c 
U(\tau')G\|_{L^{2}} d\tau'\les \tau^{\f12} \|\nab F\|_{L^{\infty}}
\|G\|_{L^{2}}
$$
which again leads to the desired estimate. 
\end{remark}
\end{proof}
\section{Sobolev space  $H^1(S)$.}
Before discussing the general, fractional,  Sobolev spaces 
in the next section it is instructive to see how the 
the standard Sobolev space $H^1(S)$ can be characterized by
 our LP   projections.
We prove  the following:
\begin{proposition} $\,$\newline {\bf i}.)\quad Consider the  LP projections
 $P_k$ associated
to an arbitrary $m\in{\cal M}_2$. Then,
\bea
\sum_k \|P_k F\|_{L^2}^2&\les& \|F\|_{L^2}^2\label{onesided1}\\
\sum_k 2^{2k}  \|P_k F\|_{L^2}^2&\les& \|\nab F\|_{L^2}^2\label{onesided2}
\eea
{\bf ii}.)\quad If in addition the LP-projections $P_k$ verify:
\be{eq:sumPk2}
\sum_kP_k^2=I
\end{equation}
Then, 
\bea
\|F\|_{L^2}^2&=&\sum_k\|P_k F\|_{L^2}^2\\
\|\nab F\|_{L^2}^2&\les&\sum_k2^{2k}\|P_k F\|_{L^2}^2
\eea
\label{prop:Sobequivalence}
\end{proposition}
\begin{proof}: \quad The first statement  of part i) is nothing else but the 
Bessel inequality established above. To prove the second statement of i) we
 write $P_k=\tilde P_k^2$ and make use of  
  of the $L^2$-finite band properties of the $\tilde P_k$'s,
 as well as  the $L^2$- boundedness of the operator $\sum_k P_k=\sum_k \tilde P_k^2$.
We shall also make use of the  following simple formula 
based on the standard definition of  
$(-\lap)^\f12$,   
\beaa
\|\nab G\|_{L^2}=<\nab G,\nab G>=-<\lap G, G>=<(-\lap)^\f12 G, (-\lap)^\f12 G>=
\|(-\lap)^\f12 G\|_{L^2}
\eeaa
Therefore,
\begin{align*}
\sum_k 2^{2k}\|P_k F\|_{L^2}^2&\les \sum_k \|\nab \tilde P_k F\|_{L^2}^2=
\sum_k \|(-\lap)^\f12
\tilde P_k F\|_{L^2}^2 = \sum_k \|
\tilde P_k (-\lap)^\f12 F\|_{L^2}^2  \\ &=
\sum_k < \tilde P_k (-\lap)^\f12 f\,, \, \tilde P_k(-\lap)^\f12 f>= 
  < \sum_k\tilde P_k^2 (-\lap)^\f12 f\,, \, (-\lap)^\f12 f> \\ &
\les\| \sum_k\tilde P_k^2 (-\lap)^\f12 F\|_{L^2} \| (-\lap)^\f12 F\|_{L^2}\les      
\|(-\lap)^\f12 F\|_{L^2}^2=\|\nab F\|_{L^2}^2
\end{align*}
as desired.
 
  The first identity of   part ii) is trivial,
$$\|F\|_{L^2}^2= < \sum_k P_k^2f, f>= \sum_k \|P_k F\|^2_{L^2}
$$
To prove the second inequality of part ii) we introduce  $P_k=\tilde P_k^2$
and  and make use of $\,\sum_k P_k^2=I$,  the $L^2$-finite band inequality
   $\|\lap \tilde P_k g\|_{L^2}\les
2^{2k}\|g\|_{L^2}$, the inequality \eqref{onesided2}, as well as 
as the commutation properties of our LP projections with $\lap$:
\beaa
\|\nab F\|_{L^2}^2&=&<-\lap f\,,\, f>=<-\lap\big(\sum_kP_k^2\,\big)  f\,,\, f>
=\sum_k<-\lap  P_kf\,,\, P_kf>\\
&\les& \sum_k\|\lap  \tilde P_k^2F\|_{L^2}\|P_k F\|_{L^2}\le 
\sum_k2^{2k}\|\tilde P_k F\|_{L^2}\|P_k F\|_{L^2}\\
&\le& \big(  \sum_k2^{2k}\|\tilde
P_k F\|_{L^2}^2 \big)^\f12\big(  \sum_k2^{2k}\|P_kF\|_{L^2}^2 \big)^\f12\les
\|\nab F\|_{L^2}\big(  \sum_k2^{2k}\|P_kF\|_{L^2}^2 \big)^\f12
\eeaa
whence,
$$\|\nab F\|_{L^2}\les \big(  \sum_k2^{2k}\|P_kF\|_{L^2}^2 \big)^\f12$$
as desired.
\end{proof}

\section{Fractional powers of $\lap$ and Sobolev spaces.}
We recall the  definition  of the Gamma function,  for  $\Re( z)>0$
\be{eq:Gamma}
\Ga(z)=\int_0^\infty e^{-t} t^{z-1} dt 
\end{equation}
as well as the beta function,
\be{eq:Beta}
B(a,b)=\int_0^1s^{a-1}(1-s)^{b-1} ds
\end{equation}
Recall that 
\be{eq:gama-beta}
B(a,b)=\frac{\Ga(a)\c\Ga(b)}{\Ga(a+b)}
\end{equation}

Let $j_a(\la)$,  $\Re(a)<0$, denote the function which is 
identically zero for $\la<0$ and 
\be{eq:jal}
j_a(\la)=\frac{1}{\Ga(-a)}\la^{-a-1}, \qquad \la>0.
\end{equation}
The following proposition is well known, 
\begin{proposition}
For all $a,b$,  $\Re(a), \,\Re(b)<0$,
\beaa
j_a*j_b=j_{a+b}
\eeaa
Moreover there exists a family of distribution $j_a$, defined for all 
$a\in {\Bbb C}$, such that,
$j_{a}*j_b=j_{a+b}$
 and $j_0=\de_0$, the Dirac delta function at the origin.
\end{proposition}
\begin{proof}:\quad
We only recall the formula $j_a*j_b=j_{a+b}$ for $\Re(a), \,\Re(b)<0$
\beaa
j_a*j_b(\la)&=&\frac{1}{\Ga(-a)}\frac{1}{\Ga(-b)}
\int_0^\la\mu^{-a-1}(\la-\mu)^{-b-1} d\mu\\
&=&\frac{1}{\Ga(-a)}\frac{1}{\Ga(-b)}\la^{-a-b-1}\
\int_0^1 s^{-a-1}(1-s)^{-b-1}\\
&=&\frac{B(-a,-b)}{\Ga(-a)\c\Ga(-b)}\la^{-a-b-1}=
\frac{1}{\Ga(-a-b)}\la^{-a-b-1}=j_{a+b}(\la)
\eeaa

\end{proof}
\begin{definition}
We define the negative  fractional powers of $\La^2=I-\lap$ on any 
 smooth tensorfield $F$ on $S$ 
according
to the formula
\be{eq:defineLaa}
\La^{a}F=\frac{1}{\Ga(-a/2)}\int_0^\infty \tau^{-\frac{a}{2}-1}
e^{-\tau}U(\tau)F d\tau
\end{equation}
 where $a$ is an arbitrary complex number with $\Re(a)< 0$.
\end{definition}
\begin{proposition}
The operators $\La^a $ is symmetric and verify the group property,
$$\La^{a}\c\La^{b} =\La^{(a+b)}$$
\end{proposition}
\begin{proof}:\quad According to the definition of $\La^a$
and the semigroup properties of $U$  we have, for any tensorfield $F$,
\beaa
\La^a\c\La^b F&=&\frac{1}{\Ga(-a/2)} \frac{1}{\Ga(-b/2)} \int_0^\infty\int_0^\infty
\tau_1^{-a/2-1}\tau_2^{-b/2-1} U(\tau_1+\tau_2)e^{-\tau_1-\tau_2} F d\tau_1 d\tau_2\\
&=& \frac{1}{\Ga(-a/2)} \frac{1}{\Ga(-b/2)}\int_0^\infty e^{-\la} U(\la) F
\int_0^\la\tau_1^{-a/2-1}(\la-\tau_1)^{-b/2-1} d\tau_1\\
 &=&\int_0^\infty j_{a/2} *j_{b/2}(\la) e^{-\la} U(\la) F=\int_0^\infty j_{a/2+b/2}(\la) e^{-\la} U(\la) F\\ &=& 
 \La^{a+b} F
\eeaa
 as desired.
\end{proof}
 We extend the definition of fractional powers of $\La$ to the range of 
 $a$ with $\Re(a)>0$, on smooth tensorfields $F$, by defining first 
 $$
 \La^a F = \La^{a-2} \c (I-\Delta) F 
 $$
 for $0< \Re (a) \le 2$ and then, in general, for 
 $0< \Re(a) \le 2m$, with an arbitrary positive integer $m$, according 
 to the formula
 $$
 \La^a F = \La^{a-2m} \c (I-\Delta)^m F.
 $$
 Observe that for $0< \Re(a)<2$,
 \beaa
\La^a F = \La^{a-2}(I-\lap)F&=&\frac{1}{\Ga(-a/2+1)}\int_0^\infty\tau^{-a/2} e^{-\tau}
 U(\tau) (I-\lap) Fd\tau\\ &=&\frac{1}{\Ga(-a/2+1)}\int_0^\infty(\frac{d}{d\tau}\tau^{-a/2}) e^{-\tau}
 U(\tau)  Fd\tau\\
&=&\frac{1}{\Ga(-a/2)}\int_0^\infty\tau^{-a/2-1} e^{-\tau}
 U(\tau)  Fd\tau
\eeaa
 Moreover, for $a=0$, the integration by parts we have performed above
yields also a boundary term. 
\beaa
\La^0 F= \La^{-2}(I-\lap)F&=&\frac{1}{\Ga(1)}\int_0^\infty e^{-\tau}
 U(\tau) (I-\lap) Fd\tau=-\int_0^\infty e^{-\tau} (I- \frac{d}{d\tau})
 U(\tau)  Fd\tau=F\\
\eeaa
i.e. $\La^{-2}$ is truly the inverse of $I-\lap$.
 \begin{remark}
 In a similar fashion, we can introduce the family of operators $\D^a= (-\Delta)^{\frac a2}$ for all $a\in {\Bbb C}$.
 As before, we start by defining formally, for $\Re (a) <0$, 
\be{eq:definecalDa}
\D^{a}F=\frac{1}{\Ga(-a/2)}\int_0^\infty \tau^{-\frac{a}{2}-1}U(\tau)F d\tau.
\end{equation}
 However, unlike $\La^a$, this formula makes sense only for smooth tensors $F$ which verify the additional 
 property that $F$ is orthogonal to the kernel of the tensor laplacean $\Delta$. 
 In view of our smoothness assumption on the manifold $S$ and the ellipticity of $\Delta$, the above
 kernel is finite dimensional.
 We can also extend the definition of $\D^a$ to the range of $a\in {\Bbb C}$ with $\Re(a) > 0$ according to 
 $$
 \D^a = \D^{a-2m} (-\Delta)^m
 $$
 with an integer $m$ such that $2m-2< \Re (a)  < 2m $.  It follows that the operators $\D^a$ with $\Re (a) >0$ 
 can naturally be extended to the space of all smooth tensors.  We also check, as before, that 
 $\D^0= I$.
 \end{remark}
  We have thus proved the following:
\begin{theorem}
There exist two family of operators   $(\D^a)_{a\in{\Bbb C}}\, ,(\La^a)_{a\in{\Bbb C}}$
such that,
\begin{align*}
&\La^0=I,\qquad   \La^{a}\c\La^{b} =\La^{(a+b)}, \qquad \La^{2k}=(I-\lap)^k,\qquad k=0,1,2\ldots,\\
&\D^0=I,\qquad \D^{2k}=(-\lap)^k, \qquad k=0,1,2\ldots
\end{align*}
on the space of all smooth tensorfields. 
In addition, the identity 
$$
\D^{a}\c\D^{b} =\D^{(a+b)}
$$ 
holds on the space of all smooth tensorfields provided that $\Re (b) > 0$ and $\Re (a+b) >0$. 
For the remaining values of $a, b\in {\Bbb C}$ the above identity holds only on the orthogonal complement of the kernel of $\Delta$.

For $\Re(a)<-2$, and any tensorfield $F$, $\La^a F$
can be defined by the formula \eqref{eq:defineLaa}, while 
$\D^a F$ is defined in  \eqref{eq:definecalDa} for $F$ in the orthogonal complement of the kernel of $\Delta$.
\label{prop:def-powers-oper}
\end{theorem}
\begin{proposition} The following estimates hold true, for any $\Re(a)<0$.
\beaa
\|\La^a F\|_{L^2}&\les &\|F\|_{L^2}.
\eeaa 
Moreover, for $2k\le \Re(a) < 2k + 2$, $k\in {\Bbb N}$,
\be{eq:funny-form-Laa}
\|(\La^a-\D^a-c_1\D^{a-2}-c_2\D^{a-4}-\cdots -c_k D^{a-2k})
F\|_{L^2}
\les \|F\|_{L^2}
\end{equation}
where $c_i=(-1)^i\frac{1}{i!}\frac{\Ga(-a/2+i)}{\Ga(-a/2)}$.
\end{proposition}
\begin{proof}:\quad
To show the boundedness of $\La^a$, $\Re(a)<0$, we only have to
use the $L^2$  boundedness of the  heat flow, 
$\|U(\tau)F\|_{L^2(S)}\le \|F\|_{L^2(S)}$.
Thus,
\beaa
\|F\|_{L^2}^{-1}\c\|\La^{-a} F\|_{L^2}&\le &       \frac{1}{\Ga(-a/2)}\int_0^\infty 
\tau^{-a/2-1} e^{-\tau} d\tau\le C_a
\eeaa
To prove \eqref{eq:funny-form-Laa} we expand $e^{-\tau}$ in the formula defining 
$\La^aF$,
$$e^{-\tau}=1 -\tau+\frac{1}{2!}\tau^2+\cdots+ (-1)^k\frac{1}{k!}\tau^k +
O(\tau^{k+1} e^{-\tau}).$$
Hence,
\beaa
\La^a F&=&\frac{1}{\Ga(-a/2)}\int_0^\infty 
\tau^{-a/2-1}e^{-\tau}U(\tau) F=\La^a F-\frac{\Ga(-a/2+1)}{\Ga(-a/2)} \La^{a-2}F\\
&+&+\frac{1}{2!}\frac{\Ga(-a/2+2)}{\Ga(-a/2)} \La^{a-4}F +\cdots
(-1)^k\frac{1}{k!}\frac{\Ga(-a/2+k)}{\Ga(-a/2)}\La^{a-2k}+E_k(F)
\eeaa
where, in view of the $L^2$ boundedness of $U(\tau)$ and the integrability
of $\tau^{-a/2+k}e^{-\tau}$ for $\Re(-a/2)+k>-1$,
we have $\|E_k(F)\|_{L^2}\les \|F\|_{L^2}$
as desired.
\end{proof}
The following proposition follows easily by standard  complex interpolation.
\begin{proposition}
 For every smooth tensorfield $F$ and any $b\ge a\ge 0$,
\bea
\|\La^a F\|_{L_2}\les\|\La^b F\|_{L^2}^{a/b}\| F\|_{L^2}^{1-a/b}\label{eq:interp-Laa}\\
\|\D^a F\|_{L_2}\les\|\D^b F\|_{L^2}^{a/b}\| F\|_{L^2}^{1-a/b}
\eea
\label{prop:interp-Laa}
\end{proposition}
We next establish  a comparison between $\|\D^a F\|_{L^2}$ and $\|\La^aF\|_{L^2}$.
\begin{proposition}
For every $a\ge 0$  and every smooth tensorfield  $F$ we have,
\beaa
\|\D^a F\|_{L^2}\les\|\La^a F\|_{L^2}\les \|\D^a F\|_{L^2}+\|F\|_{L^2}
\eeaa
\label{prop:comp-La-Da}
\end{proposition}

\begin{proof}:\quad
Indeed, according to the expansion
\eqref{eq:funny-form-Laa}, we have for $k\in {\Bbb N}$ 
for which $2k \le \Re(a)< 2k+2$,
\beaa
\|(\La^a-\D^a)F\|_{L^2}\les \sum_{i=1}^k\|\D^{a-2i}F\|_{L^2}
\eeaa
Thus, in view of the interpolation formulas of proposition \ref{prop:interp-Laa},
\beaa
\|\La^a F\|_{L^2}&\le &\|\D^a F\|_{L^2}+\sum_{i=1}^k\|\D^{a-2i}F\|_{L^2}\les
\|\D^a F\|_{L^2}+\|F\|_{L^2}
\eeaa
To prove the remaining  estimate,  
$\|\D^a F\|_{L^2}\les\|\La^a F\|_{L^2}$
it suffices to prove that, the operators
 $\La^{-a}\D^a$  are bounded in $L^2$. Observe that 
$\La^{-2}\D^2=I-\La^{-2}$. Thus, $\La^{-2}\D^2$ is bounded.
On the other hand, since
the operators $\La^a$ and $\D^a$ are selfadjoint and  commute with each other,
\beaa
\|\La^{-a}\c\DD^a F\|_{L^2}^2=<\La^{-2a}\c\DD^{2a} F, F>\les 
\|\La^{-2a}\c\DD^{2a} F\|_{L^2}\c\|F\|_{L^2}
\eeaa
Thus $\La^{-a}\D^a$ is bounded in $L^2$ if $\La^{-2a}\D^{2a}$ is.
On the other hand if  $\La^{-a}\c\D^a$,  $\La^{-b}\c\D^b$ are bounded
in $L^2$ so is $\La^{-a-b}\c\D^{a+b}$. Thus, since we already
know that $\La^{-2}\D^2$ is $L^2$ bounded, we easily infer
that $\La^{-a}\D^a$  are all bounded for all  positive  numbers 
of the form $m2^{-k}$, $m,k\in {\Bbb Z}$. The general statement 
follows now by a limiting argument.
\end{proof} We are now ready to define Sobolev norms
as follows.
\begin{definition}
For positive values of $a$ we set,
\beaa
\|F\|_{H^a(S)}=\|\La^a F\|_{L^2(S)}\approx\big(\|\D^a F\|_{L^2(S)}^2
+\|F\|_{L^2(S)}^2\big)^\f12
\eeaa
\end{definition}
In the  next theorem  we give a characterization of the Sobolev norm defined
above with the help of  LP projections. The proof depends heavily
on the following lemma:
\begin{lemma}
For all values of $a\in{\Bbb C}$  and any family of  LP projections $P_k$ with
symbol $m$  there  exists
 another family  of  LP projection $\bar{P}_k$, with symbol $\bar{m}=m*j_{a/2}$,
such that,
$$P_k \D^a F=\D^aP_k F=2^{2ak} \bar{P}_k F.$$
\label{le:nice-LP-D}
\end{lemma}

\begin{proof}:\quad 
Since the statement  is clearly true for even positive integers 
it suffices to check it  for $\Re(a)<0$. In this case,
\beaa
\D^a P_k F&=& \frac{1}{\Ga(-a/2)}\int_0^\infty \tau^{-a/2-1}U(\tau)P_k F d\tau\\
&=& \frac{1}{\Ga(-a/2)}\int_0^\infty \int_0^\infty
\tau_1^{-a/2-1}m_k(\tau_2)U(\tau_1+\tau_2)F d\tau_1 d\tau_2\\
&=&\int_0^\infty J_k(\la) U(\la) Fd\la
\eeaa
where 
\beaa
J_k(\la)&=&\frac{1}{\Ga(-a/2)}\int_0^\la m_k(\tau)(\la-\tau)^{-a/2-1}d\tau
\\&=&\frac{1}{\Ga(-a/2)} 2^{2k}\int_0^\la m(2^{2k}\tau)(\la-\tau)^{-a/2-1}d\tau\\
&=&\frac{1}{\Ga(-a/2)}\int_0^{2^{2k}\la}m(x)(\la-2^{-2k}x)^{-a/2-1}dx\\
&=&2^{ak}2^{2k}\frac{1}{\Ga(-a/2)}\int_0^{2^{2k}\la}m(x)(2^{2k}\la-x)^{-a/2-1}dx=
2^{ak}2^{2k}\bar{m}(2^{2k}\la)
\eeaa
and 
$$\bar{m}(\la)=\frac{1}{\Ga(-a/2)}\int_0^\la m(x)(\la-x)^{-a/2-1}dx=m* j_{a/2}(\la),$$
 is clearly a symbol in $\MM$.
Therefore,
$$\D P_k F=2^{2ak}\bar{P}_kF$$
as desired.
\end{proof}

\begin{theorem}$\,$\newline {\bf i}.)\quad Consider the  LP projections
 $P_k$ associated
to an arbitrary $m\in{\cal M}$. Then, for any $a\ge 0$ and any smooth
tensorfield $F$,
\bea
\sum_k 2^{2ak}  \|P_k F\|_{L^2}^2&\les& \|\D^a F\|_{L^2}^2\label{onesided2}
\eea
{\bf ii}.)\quad If in addition the LP-projections $P_k$ verify:
\be{eq:sumPk2}
\sum_kP_k^2=I
\end{equation}
then,  for\footnote{In fact the estimate holds true
 for large $a$ provided that sufficiently
many moments of the symbol $m$ of $P_k$'s vanish, see remark 
\ref{rem:inside-thm}.}
$0\le a<2$,
\bea
\|\D^a F\|_{L^2}^2&\les&\sum_k2^{2ak}\|P_k F\|_{L^2}^2
\eea
\label{thm:Sobequivalence}
\end{theorem}
\begin{proof}:\quad
 For $a=0$ part i) is nothing else but the 
Bessel inequality established earlier. To prove \eqref{onesided2}
for all $a>0$.
we make use of lemma \ref{le:nice-LP-D}. Let $P_k$
 an arbitrary family of  LP projections according  with symbol $m\in \MM$.
Let $\tilde{P}_k$  be the LP-family defined by the symbol $\bar{m}=m*j_{-a/2}$.
In view of lemma  \ref{le:nice-LP-D}\,\,
$\tilde{P}_k\D^a F=2^{ak}P_k' F$ with the corresponding symbols $\tilde{m}$
and $m'$ verifying:
 $$m'=\tilde{m}*j_{a/2}=(m*j_{-a/2})*j_{a/2}=m*(j_{-a/2}*j_{a/2})=m*\de=m.$$ 
Therefore $\tilde{P}_k\D^a F=2^{ak}P_k F$ and consequently, using Stein-Cotlar
lemma as in the proof of part iii) of theorem \ref{thm:LP},
\beaa
\sum_k\|2^{ak}P_k F\|_{L^2}^2&=&\sum_k\|\tilde{P}_k\D^a F\|_{L^2}^2=
<(\sum_k\,\tilde{P}_k^2)\,\D^{a} F\,, \,\D^a F>\\
&\le&
 \|(\sum_k\,\tilde{P}_k^2)\,\D^{a} F\|_{L^2}\|\D^a F\|_{L^2}\les\|\D^a F\|_{L^2}^2
\eeaa
as desired.

 To prove part ii)  we observe  that, if 
$$\|G\|_{L^2}^2= < \sum_k P_k^2 G, G>= \sum_k \|P_k G\|^2_{L^2}
$$
Thus, using lemma \ref{le:nice-LP-D} once more,
\beaa
\|\D^a F\|_{L^2}^2=\sum_k\|P_k\D^a F\|_{L^2}^2=
\sum_k2^{2ak}\|\bar{P}_k F\|_{L^2}^2 
\eeaa
It remains to prove that,
\be{eq:sum-Pk-barPk}\sum_k2^{2ak}\|\bar{P}_k F\|_{L^2}^2\les \sum_k2^{2ak}\|P_k
F\|_{L^2}^2
\end{equation}
To show this we proceed as follows, with the help
of the almost orthogonality estimate $\|P_{k'}\bar{P}_kG\|_{L^2}\les 2^{-2|k-k'|}
\|G\|_{L^2}$. Thus setting $J^2=\sum_k2^{2ak}\|\bar{P}_k F\|_{L^2}^2$
\beaa
J^2&=&\sum_k2^{2ak}<\bar{P}_k^2 F, F>
=\sum_{k,k'}2^{2ak}<\bar{P}_k^2 F\,,\, P_{k'}^2F>\\&=&
\sum_{k,k'}2^{2ak}<P_{k'}\bar{P}_k P_{k'}F\,, \,\bar{P}_k F>
\les\sum_{k,k'}2^{2ak}\|P_{k'}\bar{P}_k P_{k'}F\|_{L^2}\c\|\bar{P}_k F\|_{L^2}\\
&\les&\sum_{k,k'}2^{2ak}2^{-2|k-k'|}\|P_{k'}F\|_{L^2}\c\|\bar{P}_k F\|_{L^2}\\
&\les&\sum_{k,k'}2^{a(k-k')}2^{-2|k-k'|}(2^{ak'}
\|P_{k'}F\|_{L^2})\c (2^{ak}\|\bar{P}_kF\|_{L^2})\\
&\les&\big(\sum_{k'}2^{2ak'}
\|P_{k'}F\|_{L^2}^2\big)^\f12\big(\sum_{k}2^{2ak}
\|\bar{P}_{k}F\|_{L^2}^2\big)^\f12=J\c \big(\sum_{k'}2^{2ak'}
\|P_{k'}F\|_{L^2}^2\big)^\f12
\eeaa
and thus,
$$J\les \big(\sum_{k'}2^{2ak'}
\|P_{k'}F\|_{L^2}^2\big)^\f12$$
as desired.
\end{proof} 
As a corollary to theorem \ref{thm:Sobequivalence}  and proposition
\ref{prop:comp-La-Da} we derive:
\begin{corollary} For an arbitrary LP projection, $a\ge 0$ and any
smooth tensor $F$ we have,
\beaa
\sum_{k\ge 0}2^{2ak}\|P_k F\|_{L^2}^2\le \|\La^a F\|_{L^2}^2
\eeaa
Moreover, if $\sum_k P_k^2 =I$,
\beaa
\|\La^a F\|_{L^2}^2\les \sum_{k\ge 0}2^{2ak}\|P_k F\|_{L^2}^2+\|F\|^2_{L^2}
\eeaa
\label{corr:equiv-sob-Pk}
\end{corollary}
\section{ Besov spaces}
In the last section we have defined invariant Sobolev norms using
the fractional integral operators $\D^a, \La^a$ and then characterized
them with the help of the LP projections.
In this section we define invariant  Besov spaces using  directly
 the LP projections $P_k$.
\begin{definition} 
Consider the LP projections associated to
a fixed $m\in {\cal M}$ such that,
$\sum_{k} P_k^2=I$
 and  define the Besov norms, for $0\le a<2$,
 \be{eq:Besov}
\|F\|_{B^a_{p,q}}=\big(\sum_{k\ge 0} 2^{a qk}\|P_k F\|_{L^p}^q\big)^{\frac{1}{q}}+
\|F\|_{L^p}
\end{equation}
\end{definition} 
\begin{proposition}
\label{pr:SobBesov}Let the LP projections $P_k$  verify $\sum_{k} P_k^2=I$
and consider the   \eqref{eq:Besov} defined  relative to them.
Let $\tilde P_k$ any family of LP-projections  associated to an arbitrary
$\tilde m\in {\cal M}$. Then, for every $0\le a\le 1$,
\bea
\sum_{k\ge 0}2^{a k}\| \tilde P_k F\|_{L^2}&\les &\|F\|_{B^a_{2,1}}
\eea
\end{proposition}
\begin{proof}:\quad We shall use the fact that, in view of the almost orthogonality
property iii) of Theorem \ref{thm:LP} of the $P_k$'s we have 
$\|P_{k'}\tilde P_kG\|_{L^2}\les 2^{-2|k-k'|}\|G\|_{L^2}$. 
In particular, 
$$
\|P_{k} \tilde P_{<0}G\|_{L^2}\les 2^{-2k} \|G\|_{L^2},\qquad
$$
Now,
\beaa
\sum_{k\ge 0}2^{k\a}\|\tilde P_kF\|_{L^2}&\le &\sum_{k,k'\ge 0}2^{k\a}\|\tilde
P_kP_{k'}^2F\|_{L^2} + \sum_{k\ge 0}2^{k\a}\|\tilde
P_kP_{<0}^2F\|_{L^2}\\ &=&\sum_{k,k'\ge 0}2^{k\a}\|P_{k'}\tilde P_kP_{k'}F\|_{L^2} + 
\sum_{k\ge 0}2^{k\a}\|P_{<0}\tilde P_kP_{<0}F\|_{L^2}
\\
&\les& \sum_{k,k'\ge 0}2^{k\a}2^{-2|k-k'|}\|P_{k'}F\|_{L^2} + 
\sum_{k\ge 0}2^{k(\a-1)}\|P_{<0}F\|_{L^2}\\
&\les& \sum_{k'\ge 0}2^{k'\a}\|P_{k'}F\|_{L^2} + \| F\|_{L^2}=\|F\|_{B^\a_{2,1}}
\eeaa
as desired.
\end{proof}
According to corollary \ref{corr:equiv-sob-Pk} the norms $ B^a_{2,2}$
are equivalent to the Sobolev  norms $H^a$, for $0\le a<2$. For the Besov index $1$
we have the obvious inequalities,
\begin{proposition} For any smooth tensorfield $F$,
\bea
\|F\|_{B^a_{2,1}}&\les &\|F\|_{B^b_{2,1}}, \qquad a\le b\\
\|F\|_{H^a}&\les& \|F\|_{B^a_{2,1}},\qquad 0\le a\\
\|F\|_{B^a_{2,1}}&\les&\|F\|_{H^b},\qquad\,\,\,\, 0\le a<b.
\eea
\end{proposition}
\begin{proposition} The following, non-sharp, Sobolev inequality holds true
with $2<p<\infty$, $a=1-\frac 2 p$ and any tensorfield $F$,
\be{eq:LpSobolev}
\|F\|_{L^p}\les \|F\|_{B^{a}_{2,1}}
\end{equation}
\end{proposition} 
\begin{proof}:\quad 
We write $F=\sum_{k\ge 0}P_k F+P_{<0} F$. Thus, in view of the $L^p$
Bernstein inequality,
\beaa
\|F\|_{L^p}&\le &\sum_{k\ge 0}\|P_kF\|_{L^p}+\|P_{<0} F\|_{L^p}\les 
\sum_{k\ge 0}2^{k(1-\frac 2 p)}\|P_kF\|_{L^2}+\|F\|_{L^2}
\les\|F\|_{B^{\a}_{2,1}}
\eeaa
\end{proof}
\section {LP - decompositions and product estimates }

Let $P_k$ the geometric LP projections associate to an $m\in\Ml$.
We also assume that $\sum_kP_k=I$. Given a tensorfield $F$
we write, for a given $k\in \Bbb Z$
\be{eq:LPdecomp}
F=P_{<k}F+P_{\ge k} F
\end{equation}
where $P_{<k}=\sum_{l<k}P_l$, $P_{\ge k}=\sum_{l\ge k}P_l$.
Given two tensors $F,g$ and $F\c g$ some geometric product between them
we decompose,
\beaa
F\c g&=& P_{\ge k}F\c  P_{\ge k}\,G
+P_{<k}F\c P_{<k}\,G+P_{<k}F\c  P_{\ge k}\,G+P_{\ge k}F \c P_{<k}\,G
\eeaa
Thus,
\bea
P_k(F\c G)&=&\pi_k(F,\,G\,) +\si_k(F,\,G\,)+\rho_k(F,\,G\,)\label{eq:newtrichotomy}\\
\pi_k(F,\,G\,)&=&   P_k\big( P_{\ge k}F\c  P_{\ge k}G\,\big)\nn \\ 
\si_k(F,\,G\,) &=&P_k\big(P_{<k}F\c P_{<k}G\,\big)\nn\\
\rho_k(F,\,G\,)&=&P_k\big(P_{<k}F\c 
P_{\ge k}G\,\big)+P_k\big(P_{\ge k}F \c P_{<k}G\,\big)\nn
\eea
Observe that for the classical LP theory, based
on the Fourier transform, the terms $\si_k$ and $\rho_k$
are absent. Unfortunately this is not the case for our definition of 
geometric LP-projections. We shall see however that the presence of such terms
does not in any way affect the main  results that can be obtained by
the standard LP-theory. In what follows we shall apply the decomposition 
\eqref{eq:newtrichotomy} to prove a geometric version of the classical
Sobolev and Besov norm multiplication estimates.
We start with the following  

\begin{lemma}
Let $F, G \in H^{1}$  and consider \eqref{eq:newtrichotomy}.
Then, the high-high interaction term  $\pi_k(F,\,G\,)$ verifies,

$$
\sum_{k\ge 0} 2^{k} \|\pi_{k}\|_{L^{2}}
\les \|F\|_{H^1}\|G\|_{H^1}
$$
\label{le:Trich}
\end{lemma}
\begin{proof}:\quad
For $k\ge 0$ we  write, $\pi_k=\pi_k^1+\pi_k^2$  where,
$$
\pi_{k}^{1}=\sum_{k<m'<m} P_{k} \big (P_{m'} F\c P_{m} G\,\big),\qquad
\pi_{k}^{2}= \sum_{k<m'<m} P_{k} \big (P_{m} F\c P_{m'} G\,\big)
$$
By symmetry it suffices to estimate $\pi_{k}^{1}$.
Using  first the dual weak Bernstein inequality for some sufficiently large 
$p<\infty$, followed
by  Cauchy -Schwartz and then again the 
direct weak Bernstein, we obtain for any $k\ge 0$, with 
$q^{-1}+ 2^{-1}={p'}^{-1}=1-p^{-1}$,
\beaa
\|\pi_{k}^{1}\|_{L^{2}}&\les& \sum_{k<m'<m} \|P_{k}\big(P_{m'} F\c P_{m}
G\,\big)\|_{L^{2}}
\\&\les& 
2^{\frac {2k}q} \sum_{k<m'<m} \| P_{m'} F \c P_{m} G\|_{L^{p'}}\\
 &\les &
\sum_{k<m'<m} 2^{-m-m'}2^{\frac{2k}q} 
\|2^{m'} P_{m'} F\|_{L^{q}} \|2^m P_{m}G\|_{L^{2}}\\
&\les &
\sum_{k<m'<m} 2^{-m}2^{\frac{2k}q} 
2^{-m'\frac{2}{q}}\|2^{m'} P_{m'} F\|_{L^{2}} \|2^m P_{m}G\|_{L^{2}}
\eeaa
Thus, in view of the  proposition \ref{prop:Sobequivalence}
\begin{align*}
\sum_{k\ge 0} 2^{k} \|\pi_{k}^{1}\|_{L^{2}}&\les 
\sum_k\sum_{k<m'<m} 2^{-m+m'}2^{(k-m')(1+\frac 2 q)} 
\|2^{m'} P_{m'} F\|_{L^{2}} \|2^m P_{m}G\|_{L^{2}}
\\
&\les \sum_{0\le m'<m}  2^{-|m-m'|}
\|2^{m'} P_{m'} F\|_{L^{2}} \|2^m P_{m}G\|_{L^{2}}\\
&\les \|F\|_{H^1}\|G\|_{H^1}
\end{align*}
\end{proof}
We are now ready to  prove 
the following product estimates.

\begin{proposition}
\label{prop:Product1}
Let $\a,\a', \b,\b' \in (0,1)$ such that $\a + \b =\a'+\b'=1$.
Then for all tensorfields $f, g$ and 
 any $0\le \ga <1$,
\begin{equation}
\label{eq:prod'}
\|F\c G\|_{B^\ga_{2,1}}\les 
\|\La^{\a+\ga} F\|_{L^{2}} \|\La^{\b} G\|_{L^{2}}+\|\La^{\a'} F\|_{L^{2}} 
\|\La^{\b'+\ga}
G\|_{L^{2}}
\end{equation}
\end{proposition}
\begin{proof}:\quad Observe that the low frequency part
$\|P_{<0}(F\c G\,)\|_{L^2}$ can be trivially estimated
in view of the dual version of  the weak Bernstein inequality with
$q^{-1} + 2^{-1}= {p'}^{-1}=1-p^{-1}$ for some sufficiently large $p$,
$$\|P_{<0}(F\c G\,)\|_{L^2}\les \|F\c G\|_{L^{p'}}\les \|F\|_{L^q}\|G\|_{L^2}$$
followed by the Sobolev embedding \eqref{eq:LpSobolev} with $\a >\frac 2p$,
$$
\|F\|_{L^{q}}\les \|F\|_{B^{\frac 2p}_{2,1}}\les \|\La^{\a}F\|_{L^{2}}
$$
Consider now the high frequency part
$\sum_{k\ge 0} \| P_{k} (F\c G\,)\|_{L^{2}}$. Decomposing
as in \eqref{eq:newtrichotomy} we write 
$$P_{k} (F\c g)=\pi_k(F,\,G\,) +\si_k(F,\,G\,)+\rho_k(F,\,G\,)$$
The estimates for the high-high interaction  term $\pi_{k}=\pi^1_k+\pi^2_k$, $k\ge 0$ are as
follows: For $k\ge 0$ we  write, $\pi_k=\pi_k^1+\pi_k^2$;
by symmetry it suffices to estimate $\pi_{k}^{1}$.
Using  first the dual weak Bernstein inequality for some sufficiently large 
$p<\infty$, followed
by  Cauchy -Schwartz and the the direct $L^p$
Bernstein,  we obtain for any $k\ge 0$, with 
$q^{-1}+ 2^{-1}={p'}^{-1}=1-p^{-1}$,
\begin{align*}
\|\pi_{k}^{1}\|_{L^{2}}&\les \sum_{k<m'<m} \|P_{k}\big (P_{m'} F\c P_{m}
G\,\big)\|_{L^{2}}\\
&\les  2^{\frac {2k}q} \sum_{k<m'<m} \| P_{m'} F \c P_{m} G\|_{L^{p'}}\\
&\les 
\sum_{k<m'<m} 2^{-\a m'-\b m}2^{\frac{2k}q} 
\|2^{m'\a} P_{m'} F\|_{L^{q}} \|2^{m\b} P_{m}G\|_{L^{2}}\\
&\les\sum_{k<m'<m} 2^{-\a m'-\b m}2^{\frac{2k}q} 
 2^{m'(1-\frac{2}{q})}\|2^{m'\a} P_{m'} F\|_{L^{2}} \|2^{m\b} P_{m}G\|_{L^{2}}\\
&\les\sum_{k<m'<m} 2^{-\b (m-m')}2^{\frac{2(k-m')}q} 
 \|2^{m'\a} P_{m'} F\|_{L^{2}} \|2^{m\b} P_{m}G\|_{L^{2}}
\end{align*}
Thus, 
\beaa
\sum_{k\ge 0}2^{k\ga}\|\pi_{k}^{1}\|_{L^{2}}&\les&
\sum_{k<m'<m} 2^{-\b (m-m')}2^{k\ga}2^{\frac{2(k-m')}q} 
 \|2^{m'\a} P_{m'} F\|_{L^{2}} \|2^{m\b} P_{m}G\|_{L^{2}}\\
&\les&\sum_{m'<m}2^{-\b (m-m')}
\|2^{m'(\a+\ga)} P_{m'} F\|_{L^2} \|2^{m\b} P_{m}G\|_{L^{2}}
\\
&\les& \|\La^{a+\ga} F\|_{L^2}\|\La^\b G\|_{L^2}
\eeaa
since $\b>0$.

Consider now,  $\si_k(F,\,G\,)=P_{k}\big (P_{<k} F\,\c
P_{<k}G\,\big)=\si_k^1+\si_k^2$,
$$\si_k^1(F,\,G\,)=\sum_{k'< k''<k}P_{k'}F\c P_{k''}\,G,\qquad 
\si_k^2(F,\,G\,)=\sum_{k''\le k'<k}\, P_{k'}F\,\c \, P_{k''}\,G.$$
 By symmetry it suffices to estimate $\si_k^1$. 
Using the $L^2$ finite band condition
followed by  the  dual weak Bernstein inequality for $p>2$ 
sufficiently close to $p=2$ and the direct $L^{p}$ Bernstein, we estimate\footnote{We
consider only the case when the derivative affects the higher frequency; 
the other case is simpler.} with $q^{-1}+2^{-1}={p'}^{-1}$ as in the 
case of $\pi_{k}$,
\beaa
\si_k^1(F,\,G\,)&\les&\sum_{k'< k''<k} 2^{-k(1-\frac 2q)}
\|P_{k'} F\|_{L^{q}} \|\nab P_{k''} G\|_{L^{2}}\\
&\les& \sum_{k'< k''<k} 2^{k''} 2^{k'(1-\frac 2q)} 2^{-k(1-\frac 2q)}
\|P_{k'} F\|_{L^{2}} \|P_{k''} G\|_{L^{2}}\\ &\les&
\sum_{k'< k''<k}  2^{k''} 2^{(k'-k)(1-\frac 2q)} 2^{-\a k'-\b k''}
\|2^{k'\a} P_{k'} F\|_{L^{2}} \|2^{k''\b} P_{k''} G\|_{L^{2}}
\eeaa
Summing over $k$ we obtain
\begin{align*}
\sum_{k}2^{k\ga} \|\si_k(F,\,G\,)\|_{L^{2}}&\les 
\sum_{k}\sum_{k'< k''<k} 2^{k\ga} 2^{(k'-k)(1-\frac 2q)} 2^{\a (k''-k')}
\|2^{k'\a} P_{k'} F\|_{L^{2}} \|2^{k''\b} P_{k''} G\|_{L^{2}}\\
&\les \sum_{k'< k''} 2^{k''\ga} 2^{(k'-k'')(1-\frac 2q)} 2^{\a (k''-k')}
\|2^{k'\a} P_{k'} F\|_{L^{2}} \|2^{k''\b} P_{k''} G\|_{L^{2}}\\&
\sum_{k'< k''} 2^{(k'-k'')(\b -\frac 2q)} 
\|2^{k'\a} P_{k'} F\|_{L^{2}} \|2^{k''(\b+\ga)} P_{k''} G\|_{L^{2}}\\&
\les \|\La^{\a} F\|_{L^{2}} \|\La^{\b+\ga} G\|_{L^{2}}
\end{align*}
provided that $\b > \frac 2q$, which can be ensured by the choice of 
$q$, as long as $\b >0$.

We now estimate  $\rho_k(F,\,G\,)=P_k\big(P_{<k}F\c 
P_{>k}G\,\big)+P_k\big(P_{>k}F \c P_{<k}G\,\big)=\rho^1_k+\rho_k^2$.
By symmetry it suffices to estimate $\rho^1_k=
 \sum_{k'<k<m}P_k\big(P_{k'} F\c P_{m} G\,\big)$.
Arguing as in the estimate for $\si_{k}$ we use the dual weak 
Bernstein inequality followed by Cauchy-Schwartz and the 
$L^{p}$ Bernstein inequality, we obtain with 
$q^{-1}+2^{-1}={p'}^{-1}$ for a sufficiently large value of $q$,
\begin{align*}
\|\rho_k^1\|_{L^2}&\les\sum_{k'<k<m}
\|P_{k} \big (P_{k'} F\,\c\, P_{m} G\,\big)\|_{L^{2}}
\les \sum_{k'<k<m} 2^{k'} 2^{\frac {2k}q}
\|P_{k} \big (P_{k'} F\,\c\, P_{m} G\,\big)\|_{L^{p'}}\\ &\les
2^{\frac {2k}q} 2^{k'(1-\frac 2q)} \|P_{k'}F\|_{L^2} \|P_{m} G\|_{L^{2}}
\\ &\les \sum_{k'<k<m} 2^{(k-k')\frac 2q} 2^{\b(k'-m)} 
\|2^{k'\a} P_{k'}F\|_{L^2} \|2^{m\b} P_{m} G\|_{L^{2}}
\end{align*}
Now summing over $k$,
\beaa
\sum_k2^{k\ga}\|\rho_k^1\|_{L^2}&\les& \sum_k\sum_{k'<k<m} 
2^{k\ga} 2^{(k-k')\frac 2q}  
2^{\b (k'-m)}\|2^{k'\a} P_{k'}F\|_{L^2} \|2^{m\b} P_{m} G\|_{L^{2}}\\
&\les& \sum_{k'<m}2^{(\b-\frac 2q) (k'-m)} 
\|2^{k'\a} P_{k'}F\|_{L^2} \|2^{m(\b+\ga)} P_{m}
G\|_{L^{2}}\\
&\les& \|\La^\a F\|_{L^2}\|\La^{\b+\ga}G\|_{L^2}
\eeaa
 we  obtain the desired estimate
provided that $\b >\frac 2q$, which can be satisfied by the choice of
$q$, as long as $\b >0$. The corresponding estimate for $\rho^2_k$
 requires  the condition that $\a >0$.

\end{proof}

\section{The sharp Bernstein inequality}
In this section we shall prove the geometric version of the 
Bernstein inequality for arbitrary tensorfields on $M$. 
The inequality requires additional assumptions on the Gauss curvature $K$
of the manifold $M$. We shall introduce the following $L^2$- norms
depending on $K$,
\be{eq:curvLp}
K_\ga:=\|\La^{-\ga} K\|_{L^2}
\end{equation}
with $0\le \ga<1$.

\begin{theorem}
Let $S$ be a 2-d weakly regular manifold  with Gauss curvature 
$K$.
{\bf i.)}\quad For any scalar function $f$ on $S$,  $0\le \ga<1$,
 any $k\ge 0$, and an arbitrary $2\le p<\infty$,
\bea
\|P_k f\|_{L^\infty}&\les& 2^k\big(1+ 2^{-\frac kp} 
\big (K_\ga^{\frac 1{p(1-\ga)}} + K_{\ga}^{\frac 1{2p}}\big ) +
1\big)\|f\|_{L^2}\label{eq:Pkf},\label{eq:strongbernscalar}\\
\|P_{<0} f\|_{L^\infty}&\les&  
\big (1  +K_\ga^{\frac 2{p(1-\ga)}} + K_{\ga}^{\frac 1{2p}}\big)
\|f\|_{L^2}\label{eq:strong-Bern-0}
\eea

{\bf ii.)}\quad For any tensorfield  $F$ on $S$, any $k\ge 0$, 
and an arbitrary $2\le p<\infty$,
\bea
\|P_k F\|_{L^\infty}&\les & 2^k\big (1+ 2^{-\frac kp} K_0^{\frac 1p} + 
 2^{-k\frac 1{p-1}} K_0^{\frac 1{p-1}}\big ) \|F\|_{L^2},
 \label{eq:strongberntensor}\\
 \|P_{<0} F\|_{L^\infty}&\les & \big (1+ (K_0^{\frac 1p} + K_0^{\frac 1{2p}})
+  K_0^{\frac 1{p-1}}\big ) \|F\|_{L^2}.
 \label{eq:strongberntensor-0}
\eea
\label{thm:Bernstein}
\end{theorem}
\begin{proof}:\quad The proof is based on an argument
involving the product estimates 
developed in the previous section.

In view of the estimate \eqref{eq:LinftyLp}, we have for $k\ge 0$,
\bea
\|P_k F\|_{L^\infty}&\les &\|\nab^2 P_k F\|_{L^2}^{\frac 1p}
\big (\|\nab P_k F\|_{L^2}^{\frac {p-2}p} \|P_k F\|_{L^2}^{\frac 1p}  +
\|P_k F\|_{L^2}^{\frac {p-1}p}\big ) + 
\|\nab P_k F\|_{L^2}\nn\\ &\les & 
2^{k\frac {p-2}p} \|\nab^2 P_k F\|^{\frac 1p}_{L^2}\|F\|_{L^2}^{\frac {p-1}p}  +2^k\|F\|_{L^2}\label{eq:Bern1}
\eea
It remains to estimate the quantity $\|\nab^2 P_k F\|_{L^2}$.
 We  do this with the help   of  the B\"ochner identity, 
\subsection{Scalar Case} Recall that the Bochner identity for
scalars has the form,
$$
\int_S |\nab^2 g|^2 = \int_S|\lap g|^2 - \int_S K|\nab g|^2.
$$
With the help of  the  product estimates developed in the previous section
with the following choice of parameters $\a=1-\ga, \,\b=\ga$ and $\a'=\ga, \,\b'=1-\ga$,
we estimate,
\beaa
\int_S K|\nab g|^2&=&\int_S(\La^{-\ga}K)(\La^{\ga}|\nab  g|^2)\\
&\le& K_\ga\|\La^{\ga}|\nab  g|^2\|_{L^2}\les  K_\ga\|\nab  g\c \nab g \|_{B^\ga_{2,1}}\\
&\les&K_{\ga} \|\La\nab g\|_{L^2} \|\La^\ga \nab g\|_{L^2}\les K_\ga
\|\La \nab g\|_{L^2}^{1+\ga} \|\nab g\|_{L^2}^{1-\ga}
\eeaa
The last inequality follows from the condition that $\ga \le 1$
 and the interpolation inequality \eqref{eq:interp-Laa}.
Since, 
$
\|\La \nab g\|_{L^2}^2\les \int_S|\nab^2 g|^2 +\int_S |\nab g|^2
$
we  infer that
\beaa
\int_S K|\nab g|^2&\le &\f12\int_S |\nab^2 g|^2 + (K_\ga^{\frac 2{1-\ga}}+
K_\ga)
\int_S |\nab g|^2
\eeaa
Therefore,
$$
\int_S |\nab^2 g|^2\le \int_S |\lap g|^2 +\f12\int_S |\nab^2 g|^2 + 
\big (K_\ga^{\frac 2{1-\ga}} + K_{\ga}\big )\int_S |\nab g|^2
$$
This implies 
\bea
\int_S |\nab^2 g|^2&\les& \int_S |\lap g|^2  + \big (K_\ga^{\frac
2{1-\ga}} + K_{\ga}\big )\int_S |\nab g|^2\label{eq:Bochconseq}
\eea
Applying \eqref{eq:Bochconseq} to $g=P_k f$ and using the inequalities
$$\|\nab P_k f\|_{L^2} \le 2^k\|f\|_{L^2},\qquad 
\|\lap P_k f\|_{L^2} \le
2^{2k} \|f\|_{L^2}
$$
we obtain
\be{eq:secondUtau}
\|\nab^2 P_k f\|_{L^2}\les \big(2^{2k}+ 
2^k \big (K_\ga^{\frac 1{1-\ga}} + K_{\ga}^{\f12}\big ) +1
\big)\|f\|_{L^2}
\end{equation}
Combining \eqref{eq:secondUtau} with \eqref{eq:Bern1},  yields
\bea
\|P_k f\|_{L^\infty}&\les& 2^k\big(1+ 2^{-\frac kp} 
\big (K_\ga^{\frac 2{p(1-\ga)}} + K_{\ga}^{\frac 1{2p}}\big ) +
1\big)\|f\|_{L^2}\label{eq:Pkf},\\
\|P_{<0} f\|_{L^\infty}&\les&  
\big (1  +K_\ga^{\frac 2{p(1-\ga)}} + K_{\ga}^{\frac 1{2p}}\big)
\|f\|_{L^2}\label{eq:P<0f}
\eea
as desired.
\subsection{Tensor case}
We recall the B\"ochner inequality \eqref{eq:Bochner-ineq} of 
Corollary \ref{Bochner-Ineq},
$$
\|\nab^2 F\|_{L^2} \les \|\Delta F\|_{L^2} +(\|K\|_{L^2}+ \|K\|_{L^2}^\f12)\|\nab F\|_{L^2}+ \|K\|_{L^2}^{\frac p{p-1}}
 \big (\|\nab F\|_{L^2}^{\frac {p-2}{p-1}} \|F\|_{L^2}^{\frac 1{p-1}} 
+ \|F\|_{L^2}\big )
$$
Applying this to $P_k F$ we obtain
\be{eq:dyadic-Bochner}
\|\nab^2 P_k F\|_{L^2} \les \big (2^{2k} +
2^k (K_0+K_0^\f12)+  2^{k\frac {p-2}{p-1}} K_0^{\frac p{p-1}}\big )
\|F\|_{L^2}
\end{equation}
Combining \eqref{eq:dyadic-Bochner} with \eqref{eq:Bern1}
we derive
$$
\|P_k F\|_{L^\infty}\les 2^k\big (1+ 2^{-\frac kp} (K_0^{\frac 1p} + K_0^{\frac 1{2p}})
+ 2^{-k\frac 1{p-1}} K_0^{\frac 1{p-1}}  \big ) \|F\|_{L^2}
$$
as desired.
\end{proof}

\section{Sharp product estimates}
In this section we prove the sharp version of the 
product estimates of Proposition \ref{prop:Product1} involving 
Besov spaces. These estimates require an additional curvature 
assumptions which vary from the scalar to the tensor case.
The former only needs the bound on the quantity $\|\La^{-\ga} K\|_{L^2}$,
while the latter requires the finiteness of $\|K\|_{L^2}$.

Let for $0\le \ga<1$
\be{eq:constdef}
A_{\ga}:=1+ K_\ga^{\frac 1{2(1-\ga)}}
\end{equation}
denote the constants appearing in the sharp Bernstein inequalities
\eqref{eq:strongbernscalar} and \eqref{eq:strongberntensor}.
\begin{proposition}
\label{prop:Product2}
Let $S$ be a 2-d weakly regular manifold with Gauss curvature $K$.
\newline
{\bf i.)}\, For all scalar functions $f, g$,
 any $0\le \a < 2 $,  and an arbitrary $0\le \ga<1$,
\begin{align}
\|f\c g\|_{B^\a_{2,1}}&\les 
\|f\|_{B^{\a}_{2,1}} \bigg (\|g\|_{B^{1}_{2,1}}+ A_{\ga}\|g\|_{B^{
\frac 12}_{2,1}}\bigg ) \nn\\ & + 
\|g\|_{B^{\a}_{2,1}}\bigg (\|f\|_{B^{1}_{2,1}} + A_{\ga}\|f\|_{B^{
\frac 12}_{2,1}} \bigg ) \label{eq:productsc} 
\end{align}
{\bf ii.)}\, For all tensorfields $F, G$, 
 any $0\le \a <2 $, and an arbitrary $2\le p<\infty$,
\begin{align}
\|F\c G\|_{B^\a_{2,1}}&\les 
\|F\|_{B^{\a}_{2,1}} \bigg (\|G\|_{B^{1}_{2,1}} + 
A_{0}^{\frac p{p-1}} \|G\|_{B^{\f12}_{2,1}}\bigg )\nn \\ & +
\|G\|_{B^{\a}_{2,1}}\bigg (\|F\|_{B^{1}_{2,1}}  + A_{0}^{\frac p{p-1}}\|F\|_{B^{\f12}_{2,1}}\bigg )  \label{eq:productten}
\end{align}
\end{proposition}
\begin{proof}:\quad 
The proof relies on the application of the sharp Bernstein inequalities proved
in the previous section. We shall only give the arguments for the 
scalar inequality
\eqref{eq:productsc}. The modifications leading to the tensor 
inequality \eqref{eq:productten} will be obvious and follow by 
replacing the scalar Bernstein inequality \eqref{eq:strongbernscalar} 
with its tensorial version \eqref{eq:strongberntensor}.
 
As in the proof of Proposition \ref{prop:Product1} 
the low frequency part $\|P_{<0}(f\c g)\|_{L^2}$ can be trivially estimated
by means of the weak Bernstein inequality. 

Consider now the high frequency part
$\sum_{k\ge 0} \| P_{k} (f\c g)\|_{L^{2}}$. Decomposing
as in \eqref{eq:newtrichotomy} we write 
$$P_{k} (f\c g)=\pi_k(f,g) +\si_k(f,g)+\rho_k(f,g)$$
The estimates for the high-high interaction  term $\pi_{k}=\pi^1_k+\pi^2_k$, $k\ge 0$ are as
follows: For $k\ge 0$ we  write, $\pi_k=\pi_k^1+\pi_k^2$;
by symmetry it suffices to estimate $\pi_{k}^{1}$.
Using  first the dual weak Bernstein inequality for some sufficiently large 
$p<\infty$, followed
by  Cauchy -Schwartz and the the direct $L^p$
Bernstein,  we obtain for any $k\ge 0$, with 
$q^{-1}+ 2^{-1}={p'}^{-1}=1-p^{-1}$,
\begin{align*}
\|\pi_{k}^{1}\|_{L^{2}}&\les \sum_{k<m'<m} \|P_{k}\big (P_{m'} f\c P_{m}
g\big)\|_{L^{2}}\\
&\les  2^{\frac {2k}q} \sum_{k<m'<m} \| P_{m'} f \c P_{m} g\|_{L^{p'}}\\
&\les 
\sum_{k<m'<m} 2^{\frac{2k}q} 
\| P_{m'} f\|_{L^{q}} \|P_{m}g\|_{L^{2}}\\
&\les\sum_{k<m'<m}2^{\frac{2k}q} 
 2^{-\frac{2m'}{q})}\|2^{m'} P_{m'} f\|_{L^{2}} \|P_{m}g\|_{L^{2}}\\
&\les\sum_{k<m'<m} 2^{-\frac{2}q (m'-k)} 
 \|2^{m'} P_{m'} f\|_{L^{2}} \| P_{m}g\|_{L^{2}}
\end{align*}
Thus, 
\beaa
\sum_{k\ge 0}2^{k\a}\|\pi_{k}^{1}\|_{L^{2}}&\les&
\sum_{k<m'<m} 2^{-\frac 2q (m'-k)}2^{-\a(m-k)} 
 \|2^{m'} P_{m'} f\|_{L^{2}} \|2^{m\a} P_{m}g\|_{L^{2}}\\
&\les&\sum_{m'<m}2^{-\a (m-m')}
\|2^{m'} P_{m'} f\|_{L^2} \|2^{m\a} P_{m}g\|_{L^{2}}
\\
&\les& \|f\|_{B^1_{2,1}}\|g\|_{B^{\a}_{2,1}}
\eeaa

Consider now,  $\si_k(f,g)=P_{k}\big (P_{<k} f\c P_{<k}g\big)=\si_k^1+\si_k^2$,
$$\si_k^1(f,g)=\sum_{k'< k''<k}P_{k'}f\c P_{k''}g,\qquad 
\si_k^2(f,g)=\sum_{k''\le k'<k}P_{k'}f\c P_{k''}g.$$
 By symmetry it suffices to estimate $\si_k^1$. 
Using the $L^2$ finite band condition according to which 
$\|\si^{1}_{k}(f,g)\|_{L^{2}}\les 2^{-2k}\|\lap \si^{1}_{k}(f,g)\|_{L^{2}}$
we decompose 
\begin{align*}
\lap \si^{1}_{k}(f,g)&=P_{k} \sum_{k'< k''<k}\bigg  (P_{k'} f\c \lap P_{k''}g
+ 2 \nab P_{k'}f\c \nab P_{k''}g  + \lap P_{k'} f\c  P_{k''}g\bigg ) 
\\ &=\si_k^{11}(f,g) +\si_k^{12}(f,g) +\si_{k}^{13}(f,g)
\end{align*}
By symmetry it suffices to estimate the terms $\si_{k}^{11}, 
\si_{k}^{12}$.
Using the Bernstein inequality we have
\beaa
\si_k^{11}(f,g)&\les&\sum_{k'< k''<k} 2^{-2k}
\|P_{k'} f\|_{L^{\infty}} \|\lap P_{k''} g\|_{L^{2}}\\
&\les& \sum_{k'< k''<k} 2^{-2k} 2^{2k''}\big (
2^{k'} + 2^{\frac {k'}2}A_{\ga} \big )\|P_{k'} f\|_{L^{2}}  
\|P_{k''} g\|_{L^{2}}\\ &\les&
\sum_{k'< k''<k}  2^{-2k}2^{k''(2-\a)} 
\big (\|2^{k'} P_{k'} f\|_{L^{2}}+ A_{\ga}\|2^{\frac{k'}2} P_{k'} 
f\|_{L^{2}}\big ) \|2^{k''\a} P_{k''} g\|_{L^{2}}
\eeaa
Summing over $k$ we obtain for $\a <2$
\begin{align*}
\sum_{k}2^{k\a} \|\si_k^{11}(f,g)\|_{L^{2}}&\les 
\sum_{k}\sum_{k'< k''<k} 2^{-(2-\a)(k-k'')} 
\big (\|2^{k'} P_{k'} f\|_{L^{2}}+ A_{\ga}\|2^{\frac{k'}2} P_{k'} 
f\|_{L^{2}}\big ) \|2^{k''\a} P_{k''} g\|_{L^{2}}\\
&\sum_{k'< k''<k}  
\big (\|2^{k'} P_{k'} f\|_{L^{2}}+ A_{\ga}\|2^{\frac{k'}2} P_{k'} 
f\|_{L^{2}}\big ) \|2^{k''\a} P_{k''} g\|_{L^{2}}\\ & \les 
\big (\|f\|_{B^1_{2,1}} + A_{\ga}\|f\|_{B^{\frac 12}_{2,1}}\big )
\|g\|_{B^{\a}_{2,1}}
\end{align*}
To estimate $\si_{k}^{12}(f,g)$ we use the 
Gagliardo-Nirenberg inequality  \eqref{eq:GNirenberg}
$$
\|f\|_{L^{4}}\les \|\nab f\|_{L^{2}}^{\frac 12}\|f\|_{L^{2}}^{\frac 12}
+ \|f\|_{L^{2}}.
$$
Using the Gagliardo-Nirenberg estimate\footnote{We drop the low order term in the Gagliardo-Nirenberg 
inequality since we consider the case of high frequencies $k\ge 0$.} followed by the scalar Bochner 
inequality \eqref{eq:Bochconseq}
\beaa
\si_k^{12}(f,g)&\les&\sum_{k'< k''<k} 2^{-2k}
\|\nab P_{k'} f\|_{L^{4}} \|\nab P_{k''} g\|_{L^{4}}\\
&\les&\sum_{k'< k''<k} 2^{-2k}
\|\nab^{2} P_{k'} f\|_{L^{2}}^{\frac 12} \|\nab^{2} P_{k''} 
g\|_{L^{2}}^{\frac 12} \|\nab P_{k'} f\|_{L^{2}}^{\frac 12} \|\nab P_{k''} 
g\|_{L^{2}}^{\frac 12} \\
&\les& \sum_{k'< k''<k} 2^{-2k}2^{\frac {k'+k''}2}\big (
2^{ {k'}} + 2^{\frac {k'}2}A_{\ga}^{\frac 12} \big )\|P_{k'} f\|_{L^{2}}  
\big (2^{ {k''}} + 2^{\frac {k''}2}A_{\ga}^{\frac 12} \big )\|P_{k''} g\|_{L^{2}}\\ 
&\les& \sum_{k'< k''<k} 2^{-2k} 2^{\frac {3k''}2 + \frac {k''}2}\big (
2^{ {k'}} + 2^{\frac {k'}2}A_{\ga}\big )\|P_{k'} f\|_{L^{2}}  
\|P_{k''} g\|_{L^{2}}\\ 
&\les&
\sum_{k'< k''<k}  2^{-2k}2^{k''(2-\a)} 
\big (\|2^{k'} P_{k'} f\|_{L^{2}}+ A_{\ga}\|2^{\frac{k'}2} P_{k'} 
f\|_{L^{2}}\big ) \|2^{k''\a} P_{k''} g\|_{L^{2}}
\eeaa
As before, summing over $k$ we obtain for $\a <2$
\begin{align*}
\sum_{k}2^{k\a} \|\si_k^{12}(f,g)\|_{L^{2}}&\les 
\sum_{k}\sum_{k'< k''<k} 2^{-(2-\a)(k-k'')} 
\big (\|2^{k'} P_{k'} f\|_{L^{2}}+ A_{\ga}\|2^{\frac{k'}2} P_{k'} 
f\|_{L^{2}}\big ) \|2^{k''\a} P_{k''} g\|_{L^{2}}\\
&\sum_{k'< k''<k}  
\big (\|2^{k'} P_{k'} f\|_{L^{2}}+ A_{\ga}\|2^{\frac{k'}2} P_{k'} 
f\|_{L^{2}}\big ) \|2^{k''\a} P_{k''} g\|_{L^{2}}\\ & \les 
\big (\|f\|_{B^1_{2,1}} + A_{\ga}\|f\|_{B^{\frac 12}_{2,1}}\big )
\|g\|_{B^{\a}_{2,1}}
\end{align*}
and the estimate for $\si_{k}(f,g)$ follows.

We now estimate  $\rho_k(f,g)=P_k\big(P_{<k}f\c 
P_{>k}g\big)+P_k\big(P_{>k}f \c P_{<k}g\big)=\rho^1_k+\rho_k^2$.
By symmetry it suffices to estimate $\rho^1_k=
 \sum_{k'<k<m}P_k\big(P_{k'} f\c P_{m} g\big)$.
\begin{align*}
\|\rho_k^1\|_{L^2}&\les\sum_{k'<k<m}
\|P_{k'} f\|_{L^{\infty}}\|P_{m} g\|_{L^{2}}
\les \sum_{k'<k<m} \big (2^{k'} + 2^{\frac {k'}2} A_{\ga}\big ) 
\|P_{k} \big (P_{k'} f\c P_{m} g\big)\|_{L^{p'}}\\ &\les
 \|P_{k'} f\|_{L^{2}} \|P_{m} g\|_{L^{2}}
\\ &\les \sum_{k'<k<m}  2^{-m\a}
\big (\|2^{k'} P_{k'} f\|_{L^{2}} + A_{\ga}\|2^{\frac {k'}2} P_{k'} 
f\|_{L^{2}}\big ) \|2^{m\a} P_{m} g\|_{L^{2}}
\end{align*}
Now summing over $k$,
\beaa
\sum_k2^{k\a}\|\rho_k^1\|_{L^2}&\les& \sum_k\sum_{k'<k<m}  
2^{-\a(m-k)}\big (\|2^{k'} P_{k'} f\|_{L^{2}} + A_{\ga}\|2^{\frac {k'}2} P_{k'} 
f\|_{L^{2}}\big ) \|2^{m\a} P_{m} g\|_{L^{2}}\\
&\les& \big (\|f\|_{B^1_{2,1}} + A_{\ga}\|f\|_{B^{\frac 12}_{2,1}}\big )
\|g\|_{B^{\a}_{2,1}}
\eeaa
 we  obtain the desired estimate
\end{proof}
\section{Operator $\nab$ on  $B^{1}_{2,1}$ space}
Motivated by classical considerations we expect the operator
of covariant differentiation $\nab$ to act continuously in 
the scale of Besov spaces: $\,\,\nab:\, B^{s}_{2,1}\to B^{s-1}_{2,1}$
for any $s\ge 1$. The weak regularity assumptions which we impose 
on the geometry of a surface $S$ gives hope to prove this mapping 
property only for sufficiently low values of $s$. 
In this section we shall show this for the particular lowest value $s=1$.
Moreover, as in the case of the B\"ochner  and sharp Bernstein inequalities,
the regularity assumptions needed to prove
the result differ drastically dependent on whether $\nab$ is considered on the space
of scalar functions or tensorfields.  
\begin{proposition}\label{prop:nab-scal-ten}
Let $S$ be a 2-d weakly regular surface with Gauss curvature $K$ 
and let the constants $A_{\ga}$ be 
as in \eqref{eq:constdef}.

{\bf i.)} \,\, For all scalar functions $f$ and an arbitrary 
$0\le \ga<1$
\be{eq:nab-cont-scalar}
\|\nab f\|_{B^0_{2,1}}\les \|f\|_{B^1_{2,1}}+ A^2_{\ga}\|f\|_{B^0_{2,1}}.
\end{equation}

{\bf ii.)}\,\, For all tensorfields $F$ and an arbitrary $2\le p<\infty$
\be{eq:nab-cont-scalar}
\|\nab F\|_{B^0_{2,1}}\les \|f\|_{B^1_{2,1}} + 
A^{\frac {2p}{p-1}}_{0}\|f\|_{B^0_{2,1}}.
\end{equation}
\end{proposition}
\begin{proof}:\quad Once again we shall only provide the 
arguments in the scalar case. The proof of part {\bf ii.)}
is similar and relies on the tensor B\"ochner inequality
\eqref{eq:Bochner-ineq}.

We consider 
\beaa
\|\nab f\|_{B^0_{2,1}}& = &\sum_k \|P_k \nab f\|_{L^2}
\les \sum_\ell \sum_k \|P_k \nab P_\ell f\|_{L^2}\\ &= &
\sum_{\ell} \sum_{k\le \ell} \|P_k \nab P_\ell f\|_{L^2} +
\sum_{\ell} \sum_{k>\ell} \|P_k \nab P_\ell f\|_{L^2}
\eeaa
Using the dual finite band property we obtain
\beaa
\sum_{\ell} \sum_{k\le \ell} \|P_k \nab P_\ell f\|_{L^2} &\les &
\sum_{\ell} \sum_{k\le \ell} 2^k \| P_\ell f\|_{L^2} \\ &\les & 
\sum_{\ell} 2^\ell  \| P_\ell f\|_{L^2}  \sum_{k\le \ell}  2^{k-\ell}\les
\|f\|_{B^1_{2,1}}
\eeaa
It remains to estimate $\sum_{\ell} \sum_{k>\ell} \|P_k \nab P_\ell f\|_{L^2}$.
Applying the finite band property followed by the scalar B\"ochner inequality 
\eqref{eq:secondUtau}
we derive
\beaa
\|P_k \nab P_\ell f\|_{L^2}&\les & 2^{-k} \| \nab^2 P_\ell f\|_{L^2}
\les  2^{-k} ( 2^{\ell} + A^2_\ga 2^\ell ) \|P_\ell f\|_{L^2}.
\eeaa
Summing we infer that 
\beaa
\sum_{\ell} \sum_{k>\ell} \|P_k \nab P_\ell f\|_{L^2}&\les & 
\sum_{\ell} (2^\ell +A_\ga) 
\|P_\ell f\|_{L^2}\sum_{k>\ell} 2^{\ell-k}\\
&\les & \|f\|_{B^1_{2,1}} + A_\ga^2 \|f||_{B^0_{2,1}}
\eeaa
as desired.
\end{proof}

\end{document}